\documentclass[11pt]{article}

\usepackage[margin=1in]{geometry}  
\usepackage[english]{babel}
\usepackage[utf8x]{inputenc}

\usepackage[T1]{fontenc}
\usepackage{amsmath,amsthm,mathtools}
\usepackage{newpxtext,newpxmath}
\usepackage{microtype}
\usepackage{parskip}
\usepackage{natbib}
\usepackage[
  colorlinks,
  citecolor=blue,
  linkcolor=red,
  anchorcolor=red,
  urlcolor=blue
]{hyperref}
\usepackage{enumitem}
\usepackage{graphicx}
\usepackage{subcaption}
\usepackage{cleveref}

\usepackage[dvipsnames]{xcolor}
\usepackage{tcolorbox}

\newlist{settings}{enumerate}{1}
\setlist[settings]{label=\textbf{\arabic*.}, leftmargin=*}
\crefname{settingsi}{Setting}{Settings}   
\Crefname{settingsi}{Setting}{Settings}

\usepackage{mdframed}
\usepackage{kantlipsum}

\mdfdefinestyle{myenvs}{%
  hidealllines=true,%
  nobreak=true, 
  leftmargin=0pt,
  rightmargin=0pt,
  innerleftmargin=0pt,
  innerrightmargin=0pt,
}

\theoremstyle{plain}
\newmdtheoremenv[style=myenvs]{theorem}{Theorem}[section]
\newmdtheoremenv[style=myenvs]{proposition}{Proposition}[section]
\newmdtheoremenv[style=myenvs]{lemma}{Lemma}[section]
\newmdtheoremenv[style=myenvs]{example}{Example}[section]
\newmdtheoremenv[style=myenvs]{corollary}{Corollary}[section]
\theoremstyle{definition}
\newmdtheoremenv[style=myenvs]{definition}{Definition}[section]
\newmdtheoremenv[style=myenvs]{assumption}{Assumption}[section]

\theoremstyle{remark}

\newcommand{\ind}{\mathbf{1}}

\newcommand{\bern}{\textsf{Bern}}

\usepackage{xspace}

\def\##1\#{\begin{align}#1\end{align}}
\def\$#1\${\begin{align*}#1\end{align*}}

\definecolor{myblue}{rgb}{.8, .8, 1}
\definecolor{mathblue}{rgb}{0.2472, 0.24, 0.6} 
\definecolor{mathred}{rgb}{0.6, 0.24, 0.442893}
\definecolor{mathyellow}{rgb}{0.6, 0.547014, 0.24}

\newcommand{\calA}{{\mathcal{A}}}

\newcommand{\calC}{{\mathcal{C}}}

\newcommand{\calG}{{\mathcal{G}}}

\newcommand{\calO}{{\mathcal{O}}}

\newcommand{\calQ}{{\mathcal{Q}}}

\newcommand{\calX}{{\mathcal{X}}}
\newcommand{\calY}{{\mathcal{Y}}}
\newcommand{\calZ}{{\mathcal{Z}}}

  \def\EE{\mathbb{E}}

  \def\NN{\mathbb{N}}
  
  \def\PP{\mathbb{P}}
  
  \def\RR{\mathbb{R}}

\def\vareps{\varepsilon}

\def\hat{\widehat}
\def\tilde{\widetilde}

\newcommand{\aci}{\texttt{ACI}}
\newcommand{\cp}{\texttt{CP}}

\DeclareTextFontCommand{\texttt}{\ttfamily\upshape}

\newcommand{\err}{\mathsf{e}}
\DeclareRobustCommand{\bbar}{\bar{\beta}}
\DeclareRobustCommand{\bunder}{\underline{\beta}}

\newcommand{\itv}{\kappa,\beta}
\newcommand{\Witv}{\Gamma_{\itv}}

\newcommand{\Corrid}{\kappa}

\def\err{\texttt{e}}

\def\miscov{\texttt{E}}

\def\Asym{\calA_{\beta,\kappa}}
\def\Aasym{\calA_{b,d,c_{-},c_+}}

\def\inst{\Delta^{\rm inst}_T}

\title{The cost of enforcing long-run coverage: a phase transition in online conformal prediction}
\author{%
  Yu Gui\textsuperscript{1} \qquad Rina Foygel Barber\textsuperscript{2} \qquad Cong Ma\textsuperscript{2}\\[1em]
  {\small \textsuperscript{1}Department of Statistics and Data Science, The Wharton School, University of Pennsylvania}\\
  {\small \textsuperscript{2}Department of Statistics, University of Chicago}%
}
\date{\today}

\begin{document}

\maketitle

\addtocontents{toc}{\protect\setcounter{tocdepth}{-3}}

\begin{abstract}
Consider an online prediction setting, where at each time, we are given a new feature vector and asked to provide a prediction interval or prediction set that will likely contain the corresponding response. In this setting, if the data points are i.i.d.\ (or more generally, exchangeable), conformal prediction provides a mechanism for building valid prediction sets around the output of any trained model, offering coverage probability $\geq 1-\alpha$ at each time $t$. In practice, however, issues such as distribution drift, changepoints, or other challenges may cause violations of the i.i.d.\ assumption. In this adversarial setting, online conformal prediction methods can still offer guarantees on long-run empirical coverage. However, there is a trade-off between these two types of goals: enforcing tight constraints on long-run empirical coverage to protect against the adversarial setting can lead to unstable oscillations in the instantaneous training-conditional coverage in the i.i.d.\ setting. 

In this work, we seek \emph{best-of-both-worlds} algorithms that aim to minimize cumulative error in instantaneous training-conditional coverage when in the i.i.d.\ setting, while also offering long-run empirical coverage guarantees in the adversarial setting. Specifically, consider a condition that requires long-run empirical coverage to deviate from the target level $1-\alpha$ by at most $\mathcal{O}(T^{-\beta})$, for some $\beta\in[0,1]$. We prove that this family of conditions can be described with a phase transition. For $\beta\in[0,\frac{1}{2}]$, it is possible to construct algorithms satisfying such a constraint in the adversarial setting, while offering the same instantaneous coverage performance as split conformal prediction if the data is in fact i.i.d.\ On the other hand, for $\beta\in(\frac{1}{2},1]$, this type of condition comes at a cost, and does not achieve best-of-both-worlds behavior. 


\end{abstract}


\section{Introduction}
\label{sec:intro}

Consider an online uncertainty quantification problem in which each data point consists of a feature $X_t\in\calX$ and a response $Y_t\in\calY$.
For each time $t\geq 1$, let $\mathcal{C}_t$ denote any prespecified construction that maps from an observed data stream to a prediction set for the data point at time $t$.
Concretely, we write
\begin{align}\label{ass:pred}
    \calC_t(X_t) = \calC_t\big(X_t; (X_1,Y_1),\dots,(X_{t-1},Y_{t-1})\big)\subseteq\mathcal{Y}
\end{align}
to denote the prediction set at time $t$, when trained on data points $(X_1,Y_1),\dots,(X_{t-1},Y_{t-1})$.
The current outcome $Y_t$ is not available when the set is constructed.
Throughout this paper, we use the notation $\calC = (\calC_t)_{t \geq 1}$ to denote the online algorithm producing the sequence of prediction sets $(\calC_t(X_t))_{t \geq 1}$.

After $Y_t$ is revealed, the performance of the prediction set can be summarized by a binary miscoverage indicator
\[
\err_t = \ind\{Y_t \notin \calC_t(X_t; (X_1,Y_1),\dots,(X_{t-1},Y_{t-1}))\}.
\]
We expect the prediction sets to frequently cover the true response $Y_t$, quantified by the nominal coverage level $1-\alpha \in (0,1)$.
To quantify this, we may require that the probability of coverage at each time point is calibrated---that is, that $\PP(\err_t \neq 1)$ or the coverage rate conditional on the history is $\approx 1-\alpha$, at each time $t$. This type of guarantee can only be achieved with some assumptions on the data stream, however. More generally, we can instead aim for a weaker guarantee, requiring that the empirical \emph{average} coverage rate is $\approx 1-\alpha$ over a period of time.
We formally define and discuss these two guarantees in the following subsections.

\subsection{The i.i.d. setting: instantaneous training-conditional coverage}
When the data points $(X_t, Y_t)$ are independent and identically distributed (or exchangeable), split conformal prediction uses past scores to calibrate prediction sets.
The score is fixed or fitted to independent data, and higher scores indicate greater deviation from the fitted model.
Concretely, with a fixed score function $s: \calX \times \calY \rightarrow \RR$, the split conformal prediction set can be defined by
\begin{equation}\label{eq:conformal-baseline-set}
 \calC_t^{\cp}(x)=
 \left\{y \in \calY:\frac{1+\sum_{r<t}\ind\{s(X_r,Y_r)\ge s(x,y)\}}{t}>\alpha\right\}.
\end{equation}
Exchangeability gives a finite-sample marginal coverage guarantee, i.e., $\PP(Y_t \in \calC_t^{\cp}(X_t)) \geq 1 - \alpha$.
A more local requirement asks whether coverage is close to $1-\alpha$ given the observed past at each time $t$.
For i.i.d.\ observations and a nonatomic score distribution,\footnote{A distribution is nonatomic if its probability measure $\mu$ satisfies that for any measurable set $A$ with $\mu(A) > 0$, there exists a measurable subset $B$ of $A$ such that $\mu(A) > \mu(B) > 0$.} the split conformal prediction set in \eqref{eq:conformal-baseline-set} satisfies the following bound \citep{vovk2012conditional,angelopoulos2024theoretical}:
\begin{equation}\label{eq:iid-tcc}
 \EE_{\PP}\left|
 \PP\{Y_t\in\calC_t^{\cp}(X_t)\mid (X_s,Y_s)_{s < t}\}-(1-\alpha)
 \right|
 \le\frac{1}{2\sqrt{t+1}}+\frac1t.
\end{equation}
Thus, under these assumptions, the instantaneous training-conditional coverage of split conformal prediction, i.e., its coverage at time $t$ given the history $(X_s,Y_s)_{s<t}$, deviates from the target level $1-\alpha$ by $O(t^{-1/2})$ in expectation. To measure this deviation over $T$ rounds, for any online algorithm $\calC$ we define
\begin{equation}\label{eq:instant}
  \inst(\calC) := \sum_{t=1}^T \EE_{\PP}\left|\PP\{Y_t\in\calC_t(X_t)\mid (X_s,Y_s)_{s<t}\}-(1-\alpha)\right|.
\end{equation}
With $\calC = \calC^{\cp}$, summing \eqref{eq:iid-tcc} over $t$ gives $\inst = O(\sqrt{T})$.

\subsection{The adversarial setting: long-run empirical coverage}
When exchangeability is not assumed, one widely adopted quantity that measures the \emph{adversarial} performance of online algorithms is the long-run cumulative coverage deviation \citep{gibbs2021adaptive}:
\begin{equation}\label{eq:ST_intro}
  \sum_{t=1}^T \left[\ind\{Y_t \in \calC_t(X_t)\} - (1-\alpha)\right] = \sum_{t=1}^T(\alpha - \err_t).
\end{equation}
Here, we denote $\miscov_T = \sum_{t=1}^T(\err_t - \alpha)$, where $|\miscov_T| / T$ tracks how far the long-run empirical coverage rate departs from $1-\alpha$.
A well-calibrated algorithm keeps $|\miscov_T|$ small relative to $T$.
This guarantee can be achieved by a large class of algorithms without imposing distributional assumptions on the data sequence $\{(X_t,Y_t), t \geq 1\}$.
To quantify the long-run empirical coverage guarantee in the adversarial setting, we say the algorithm satisfies a \emph{$(\beta,\Corrid)$-constraint} if
\begin{equation}\label{eq:corridor_intro}
  |\miscov_T| / T \leq \Corrid T^{-\beta}
  \quad \text{for all } T \ge 1.
\end{equation}
Throughout this paper, $\beta\in[0,1]$ and $\kappa\ge1/2$, so the interval has width at least one and a feasible error outcome exists at every round.
Note that if $\beta=0$ and $\kappa \geq \max\{\alpha,1-\alpha\}$, then the constraint holds trivially for any sequence $(\err_t)_{t \geq 1}$.
For $\beta>0$, this implies that the long-run empirical coverage rate converges to $1-\alpha$.
A larger $\beta$ requires tighter control over the long-run empirical coverage rate.
We use $\Asym$ to collect all algorithms that satisfy \eqref{eq:corridor_intro}.
\footnote{For simplicity of presentation, we present the main theoretical results under symmetric constraints and give their asymmetric counterparts in Appendix~\ref{app:asymmetric}. We can likewise define asymmetric constraints: we say the algorithm satisfies a \emph{$(b,d,c_{-}, c_{+})$-constraint} if
$\miscov_T / T \;\in\; [-c_{-}T^{-d},\; c_{+}T^{-b}]$.
We use $\Aasym$ to collect such algorithms. The same phase transition holds in the asymmetric setting, with $\min\{b,d\}$ playing the role of $\beta$: the minimax oracle disagreement has order $T^{\max\{1/2,\min\{b,d\}\}}$.}

Our central question is how close an algorithm in $\Asym$ can stay to exact coverage given the history when the data are i.i.d.
By \eqref{eq:iid-tcc}, the instantaneous training-conditional coverage of split conformal prediction is within $O(t^{-1/2})$ of $1-\alpha$ in expectation at round $t$, so its cumulative deviation sum to $O(\sqrt T)$ over $T$ rounds, i.e., $\inst= O(\sqrt{T})$ for $\calC = \calC^{\cp}$.
For a nontrivial constraint, split conformal prediction, however, need not belong to $\Asym$.
We show that the smallest cumulative error in the instantaneous training-conditional coverage that an algorithm in $\Asym$ can guarantee for i.i.d.\ data is of order
\begin{equation}\label{eq:preview-rate}
 T^{\max\{1/2,\beta\}}=
 \begin{cases}
  \sqrt T, & 0\le\beta\le1/2,\\
  T^{\beta}, & 1/2<\beta\le1.
 \end{cases}
\end{equation}
Thus, the constraint incurs no cost in rate when $\beta\le1/2$, whereas for $\beta>1/2$ the rate is dominated by the constraint, yielding a phase transition at $\beta=1/2$.

\subsection{Adaptive conformal inference (ACI) revisited}

Adaptive Conformal Inference (ACI), introduced by \citet{gibbs2021adaptive}, is one way to enforce long-run empirical coverage on every data sequence.
We consider its quantile-threshold form \citep{angelopoulos2024online,angelopoulos2024theoretical}, which uses a fixed score $s:\calX\times\calY\to\RR$ and constructs $\calC_t^{\aci}(x)=\{y:s(x,y)\le q_t\}$.
After observing the outcome, ACI updates its threshold by
\begin{equation}\label{eq:aci}
  q_{t+1} \;=\; q_t
  + \eta_t\bigl(\err_t - \alpha\bigr),
\end{equation}
where $\eta_t > 0$ is the step size.
We note that the original ACI procedure \citep{gibbs2021adaptive} updates the nominal coverage level. Here we analyze the quantile-threshold update \eqref{eq:aci}.
The threshold increases after miscoverage and decreases after coverage.
The step sizes determine how tightly ACI controls the cumulative coverage deviation in the adversarial case.
More concretely, when $s(x,y)\in[0,B]$, $\eta_t\equiv\eta$ (i.e., a constant step size), and $q_1\in[0,B]$, ACI has the long-run empirical coverage guarantee in the adversarial setting \citep{gibbs2021adaptive,angelopoulos2024online}
\[
|\miscov_T| = \big| \sum_{t = 1}^T (\err_t - \alpha) \big| \le \frac{B+\eta}{\eta},
\]
indicating that ACI satisfies the constraint with $\beta = 1$ and a constant sufficiently large to dominate $(B+\eta)/\eta$.



\subsection{Best-of-both-worlds algorithms?}

Recall that split conformal prediction satisfies \eqref{eq:iid-tcc} but can violate \eqref{eq:corridor_intro} on some data sequences.
Constant-step ACI satisfies \eqref{eq:corridor_intro}, but its threshold moves by at least $\eta\min\{\alpha,1-\alpha\}$ after each observation, so its instantaneous training-conditional coverage continues to fluctuate even for i.i.d.\ data.
Since it is rarely known in advance whether the data are i.i.d., we seek a single algorithm that satisfies the long-run constraint for every data sequence and, when the data are i.i.d., keeps its instantaneous training-conditional coverage as close as possible to $1-\alpha$.
We call such an algorithm \emph{best-of-both-worlds}.

For ACI, this goal suggests preserving its long-run empirical coverage bound while reducing the size of each threshold update, that is, specifying decaying step sizes \citep{angelopoulos2024online}.
Concretely, when $s(x,y)\in[0,B]$, $q_1\in[0,B]$, and $\eta_t=ct^{-(1-\beta)}$ with $c>0$ and $0<\beta<1$, the bound of \citet[Theorem~1]{angelopoulos2024online} gives
\begin{align}\label{ex:aci_decay}
 \frac{|\miscov_T|}{T}\le\frac{B+c}{c}\,T^{-\beta}
 \qquad\text{for every }T\ge1.
\end{align}
Thus ACI satisfies a $(\beta,\kappa)$-constraint with $\kappa=(B+c)/c$ on every data sequence.
Faster decay of the step sizes corresponds to smaller $\beta$ and a weaker long-run empirical coverage guarantee.
Under i.i.d.\ sampling with a fixed bounded score whose CDF is continuous and whose $(1-\alpha)$-quantile is unique, the steps $\eta_t=ct^{-(1-\beta)}$ with $0<\beta<1/2$ make the working quantile $q_t$ converge almost surely to the population $(1-\alpha)$-quantile \citep[Theorem~3]{angelopoulos2024online}.
Instantaneous training-conditional coverage then converges to $1-\alpha$ as well.
Together with \Cref{ex:aci_decay}, this gives both guarantees asymptotically.

Each new observation changes ACI's threshold by an amount proportional to $\eta_t$, even as the threshold converges to the population quantile.
A larger $\beta$ yields tighter long-run empirical coverage but causes these updates to decay more slowly, raising the question of how much error in the instantaneous training-conditional coverage accumulates over time.
More concretely, for i.i.d.\ data under the score-density assumptions in \Cref{sec:aci-subopt}, \Cref{prop:aci-subopt} quantifies this cost.
With $\eta_t=ct^{-(1-\beta)}$ and $0<\beta<1$, it gives
\begin{align}
 \sum_{t=1}^T\EE_{\PP}\left|
 \PP\{Y_t\in\calC_t^{\aci}(X_t)\mid (X_s,Y_s)_{s < t}\}-(1-\alpha)
 \right|\asymp T^{(1+\beta)/2}.
\end{align}
For every $0<\beta<1$, this rate exceeds the $O(\sqrt {T})$ bound for split conformal prediction in \eqref{eq:iid-tcc}.
Split conformal prediction, however, does not itself enforce the long-run adversarial constraint on every sequence. Then it motivates us to ask the question: \emph{what is the smallest cumulative error in instantaneous training-conditional coverage achievable under the long-run empirical coverage constraints?}

\begin{figure}[t]
\centering
\includegraphics[height=0.35\textwidth]{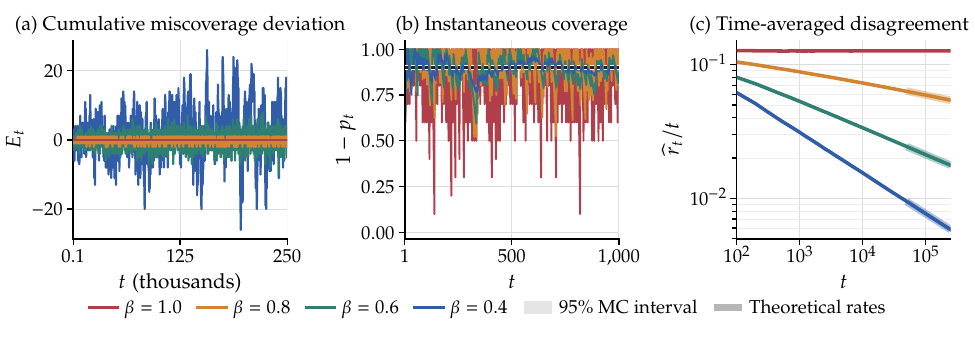}
\caption{\textbf{Long-run empirical coverage versus oracle disagreement in ACI.}
We use the regression model in \Cref{sec:simulations}, with $\alpha=0.1$, $q_1=0.9$, and step sizes $\eta_t=t^{-(1-\beta)}$.
We define $p_t=\PP(Y_t\notin\calC_t^{\aci}(X_t)\mid (X_s,Y_s)_{s < t})$ and compute $\widehat r_t$ by averaging $\sum_{s\le t}|p_s-\alpha|$.
Dashed guides are proportional to $t^{-(1-\beta)/2}$ for $\beta<1$, as predicted by \Cref{prop:aci-subopt}.}
\label{fig:aci_motivation}
\end{figure}

\Cref{fig:aci_motivation} illustrates the tradeoff in a synthetic dataset (\Cref{sec:simulations}), where the fixed score, by construction, is uniformly distributed on $[0,1]$.
We compare $\eta_t=t^{-(1-\beta)}$ for $\beta\in\{1,0.8,0.6,0.4\}$, including constant steps at $\beta=1$.
In panel (c), we focus on a metric that quantifies the expected disagreement between ACI prediction sets and the population score oracle; a formal definition is presented in \Cref{sec:main}.
In this experiment, a larger $\beta$ keeps the cumulative coverage deviation closer to zero, but instantaneous training-conditional coverage fluctuates more, and the time-averaged oracle disagreement decreases more slowly.
These observations leave open the extent to which ACI's statistical cost is required by the coverage constraint and how much arises from its own algorithmic design.

\paragraph{Core question.}
Motivated by the aforementioned tradeoff, it remains an open question how much instantaneous training-conditional coverage error is unavoidable under the long-run empirical coverage constraint.
Our goal is to find an algorithm that offers the \emph{best-of-both-worlds} guarantee: it satisfies the prescribed long-run empirical coverage constraint on every data sequence and achieves the best statistical performance when the data are i.i.d.

To make this objective precise, we compare an algorithm with an oracle that knows the data distribution and has coverage exactly $1-\alpha$ given the past.
Their cumulative disagreement counts the rounds on which the observed response belongs to exactly one of their prediction sets.
Our central question is:
\begin{quote}
How closely can a single algorithm match an exact oracle under i.i.d.\ sampling while satisfying a prescribed long-run empirical coverage constraint on every data sequence?
\end{quote}
In particular, can the $\sqrt T$ cumulative rate of split conformal prediction be retained, and when does tighter long-run empirical coverage force a higher cost?
We seek the smallest possible cost and an algorithm that attains it without knowing the data distribution or whether the observations are i.i.d.

\paragraph{Paper organization.}
\Cref{sec:main} formalizes oracle disagreement and separates the difficulty of approximating an oracle from the cost of the long-run empirical coverage constraint.
For the family of all exact oracles, under i.i.d.\ laws with a nonatomic distribution for the fixed score, we prove a minimax rate of $T^{\max\{1/2,\beta\}}$ and achieve it by truncating split conformal prediction only when needed to enforce the constraint.
Thus the $\sqrt T$ rate is retained for $\beta\le1/2$, while tighter constraints necessarily increase the rate.
\Cref{sec:aci-subopt} shows that fixed-score ACI with steps proportional to $t^{\beta-1}$ incurs a larger cost for $0<\beta<1$ under the stated density assumptions.

\section{Oracle disagreement under an adversarial constraint}
\label{sec:main}

Recall the setting introduced in \Cref{sec:intro}.
We use $\calA$ to denote the class of measurable online algorithms in \eqref{ass:pred} that are independent of the unknown data distribution.
At round $t$, an online algorithm $\calC$ uses the past observations $(X_r,Y_r)_{r<t}$ and the current feature $X_t$ to construct a prediction set $\calC_t(X_t)\subseteq\calY$ before observing $Y_t$.
Write $\err_t^{\calC}=\ind\{Y_t\notin\calC_t(X_t)\}$ and $\miscov_T^{\calC}=\sum_{t=1}^T(\err_t^{\calC}-\alpha)$.
With a fixed joint distribution $\PP$ for the data sequence $(X_t,Y_t)_{t \in \NN}$ and to make the conditioning clear, we define $\calG_{t-1} = \sigma((X_r,Y_r)_{r<t})$.\footnote{When needed, we use $\zeta_{1:t} = (\zeta_1, \dots, \zeta_t)$, which is independent of the data, to collect auxiliary randomization up to time $t$. Throughout this paper, all algorithms and oracles considered here are deterministic measurable functions of the observed data.}

In this section, we formally quantify the instantaneous training-conditional performance of an online algorithm. To define this metric, we compare each prediction set $\calC_t(x)$ at time $t$ with an \emph{oracle} whose instantaneous training-conditional coverage is exactly calibrated at each $t$.
Concretely, an oracle $\calC^*_{\PP}$ may know $\PP$, uses the same available information as $\calC$, and satisfies
\begin{equation}\label{eq:main-oracle-calibration}
\textnormal{Definition of an oracle:}\qquad
 \PP\bigl(Y_t\in\calC^*_{\PP,t}(X_t)\mid\calG_{t-1}\bigr)=1-\alpha
 \qquad\text{almost surely for every }t\ge1.
\end{equation}
As an example, under i.i.d.\ sampling with a measurable score function $s:\calX \times \calY \rightarrow \RR$, if the distribution of $s(X_1,Y_1)$ is nonatomic, let $F_{\PP,s}$ be its CDF and set
\begin{equation}\label{eq:population-oracle}
 q_{\PP,s}=\inf\{q:F_{\PP,s}(q)\ge1-\alpha\},\qquad
 \calC^*_{\PP,s,t}(x)=\{y:s(x,y)\le q_{\PP,s}\}.
\end{equation}
Independence of the data sequence and continuity of $F_{\PP,s}$ give the oracle condition \eqref{eq:main-oracle-calibration}.
We note that, in general, the existence of $\calC^*_{\PP}$ requires regularity conditions on the joint distribution $\PP$, and for simplicity, we assume $\PP$ is nonatomic throughout this paper unless otherwise specified.

To evaluate two online algorithms, we measure their cumulative disagreement on the same data sequence by
\begin{equation}\label{eq:main-loss}
 \ell_T(\calC,\calC')=\sum_{t=1}^T|\err_t^{\calC}-\err_t^{\calC'}|.
\end{equation}
This counts responses that belong to exactly one of the two sets.
Then, the expected disagreement with an oracle $\calC^*_{\PP}$, i.e., $\EE_{\PP}[\ell_T(\calC,\calC^*_{\PP})]$, quantifies the instantaneous training-conditional performance of an algorithm $\calC = (\calC_t)_{t \in \NN}$ under a given data distribution $\PP$.
For $p_t^{\calC}=\PP(\err_t^{\calC}=1\mid\calG_{t-1})$, conditional Jensen's inequality further gives
\begin{equation}\label{eq:coverage-disagreement}
 \EE_{\PP}\ell_T(\calC,\calC_{\PP}^*)
 \ge\sum_{t=1}^T\EE_{\PP}|p_t^{\calC}-\alpha| = \inst
 \qquad \textnormal{for every exact oracle }\calC_{\PP}^*.
\end{equation}
Equality holds when two prediction sets are coupled, i.e., the sets are defined by two thresholds on the same score, whose cutoffs are determined by past data.

\subsection{Hardness of matching an oracle under long-run empirical coverage constraints}
\label{sec:hardness}

Let $\tilde\calQ$ be the class of all i.i.d.\ joint laws on the data sequence $(X_t,Y_t)_{t\in\NN}$.
Fix a measurable score $s:\calX\times\calY\to\RR$, and let
\[
 \calQ=\{\PP\in\tilde\calQ:
          s(X_1,Y_1)\text{ has a nonatomic distribution under }\PP\}.
\]
We assume that $\calQ$ is nonempty.
For each $\PP\in\calQ$, let $\calO_\alpha(\PP)$ denote the set of algorithms that satisfy the oracle condition \eqref{eq:main-oracle-calibration}.
Particularly, this family includes the population score oracle defined in \eqref{eq:population-oracle}.

We measure an algorithm's performance by
\begin{equation}\label{eq:generic-risk}
 \sup_{\PP\in\calQ}\inf_{\calC^*_{\PP}\in\calO_\alpha(\PP)}
 \EE_{\PP}\ell_T(\calC,\calC^*_{\PP}).
\end{equation}
The infimum allows a separate oracle for each algorithm and law, without giving the algorithm access to the oracle or the distribution.
Throughout, $\lesssim$, $\gtrsim$, and $\asymp$ compare quantities up to positive constants for all sufficiently large $T$.

Matching an oracle is difficult for two reasons.
The algorithm must learn from an unknown distribution and keep its realized errors within the prescribed bounds.
The following result separates these two costs.

\begin{theorem}\label{thm:generic-hardness}
Under the setup above, every $\calC\in\Asym$ satisfies
\begin{equation}\label{eq:generic-lower}
\sup_{\PP\in\calQ}
 \inf_{\calC^*_{\PP}\in\calO_\alpha(\PP)}
 \EE_{\PP}\ell_T(\calC,\calC^*_{\PP}) \geq C_{\alpha,\beta,\kappa} T^{\max\{1/2,\beta\}},
\end{equation}
for sufficiently large $T$, where $C_{\alpha,\beta,\kappa} > 0$ depends only on $\alpha, \beta, \kappa$.
\end{theorem}
More generally, the lower bound can be decomposed as
\begin{equation}\label{eq:generic-lower-decomposition}
\begin{split}
 &\sup_{\PP\in\calQ}
 \inf_{\calC^*_{\PP}\in\calO_\alpha(\PP)}
 \EE_{\PP}\ell_T(\calC,\calC^*_{\PP})\\
 &\qquad\gtrsim
 \underbrace{
 \inf_{\calC^0\in\calA}\sup_{\PP\in\calQ}
 \inf_{\calC^*_{\PP}\in\calO_\alpha(\PP)}
 \EE_{\PP}\ell_T(\calC^0,\calC^*_{\PP})
 }_{\asymp\sqrt T}
 +
 \begin{cases}
  0,&0\le\beta<1/2,\\
  T^\beta,&1/2\le\beta\le1.
 \end{cases}
\end{split}
\end{equation}
The first term is the unavoidable cost of learning from past data.
Even without a long-run empirical coverage constraint, no algorithm can achieve disagreement smaller than order $\sqrt T$ uniformly over $\calQ$.
This remains true when comparison is allowed with every exact oracle.
Split conformal prediction attains this order by estimating the population score quantile.

The second term is the cost of enforcing the long-run empirical coverage constraint.
Exact calibration does not prevent the oracle's realized errors from fluctuating: its centered error sum has standard deviation $\sqrt{\alpha(1-\alpha)T}$.
However, the constraint allows a deviation of only $\kappa T^{1-\beta}$.
When $\beta\ge1/2$, keeping the errors within these bounds requires repeated departures from the oracle.
These departures force a cumulative disagreement with the order of at least $T^\beta$, even if the distribution is known.

The higher cost determines the lower-bound rate, which exhibits a sharp phase transition at $\beta=1/2$.
The cost of learning the unknown distribution dominates when $\beta<1/2$, both costs have order $\sqrt T$ at $\beta=1/2$, and the long-run empirical coverage constraint dominates when $\beta>1/2$.
In the next subsection, we construct an algorithm that achieves this rate in both regimes.

If distributions are not nonatomic, there may be no oracle that satisfies the condition \eqref{eq:main-oracle-calibration}, but we can approach the problem differently by considering relaxed oracles; see \Cref{app:relaxed}.

\subsection{Truncation attains the minimax rate}\label{sec:score-minimax}
In this section, we present an explicit online construction of the online algorithm that has a matching upper bound in \Cref{thm:generic-hardness}.
Fix a measurable score $s:\calX\times\calY\to\mathbb R$.
Recall that $\calO_\alpha(\PP)$ denotes the set of oracle algorithms that satisfy \eqref{eq:main-oracle-calibration} and under $\PP \in \calQ$, the distribution of score $s(X_t,Y_t)$ is nonatomic.

To meet the long-run empirical coverage constraints, a straightforward way is to truncate the cumulative miscoverage error $\miscov_T$ and enforce the adversarial constraints at each step.
Concretely, for any baseline $\calC^0$, define $\widehat\calC=\Witv\circ\calC^0$ recursively, starting at $\miscov_0^{\widehat\calC}=0$, by
\begin{equation}\label{eq:intervention}
 (\Witv \circ \calC^0)_t(x) = \widehat\calC_t(x)=
 \begin{cases}
  \calY,&\miscov_{t-1}^{\widehat\calC}+1-\alpha>\kappa t^{1-\beta},\\
  \varnothing,&\miscov_{t-1}^{\widehat\calC}-\alpha<-\kappa t^{1-\beta},\\
  \calC_t^0(x),&\text{otherwise}.
 \end{cases}
\end{equation}
If a potential miscoverage at round $t$ would push the cumulative deviation above its upper bound, we output $\calY$ to ensure coverage.
Similarly, if coverage would push $\miscov_t^{\widehat\calC}$ below its lower bound, we output $\varnothing$ to force a miscoverage.
Otherwise, either outcome respects the bounds, so we keep the baseline prediction set.
$\miscov_t^{\widehat\calC}$ is computed from the errors of $\widehat\calC$, whereas the baseline $\calC^0$ continues to use every observation, including those from rounds in which its set is replaced by its truncated counterpart.
For $\kappa\ge1/2$, the two forcing conditions cannot hold simultaneously, and the resulting prediction sets satisfy \eqref{eq:corridor_intro} on every data sequence.

With split conformal prediction as the baseline, this truncation attains the minimax rate of expected oracle disagreement, as the following theorem shows.
\begin{theorem}\label{thm:optimal}
The single algorithm $\Witv\circ\calC^{\cp}\in\Asym$ attains the upper bound:
\begin{equation}\label{eq:rates}
 \sup_{\PP\in \calQ}\inf_{\calC^*_{\PP}\in\calO_\alpha(\PP)}
 \EE_{\PP}\ell_T(\Witv\circ\calC^{\cp},\calC^*_{\PP})
 \lesssim T^{\max\{1/2,\beta\}}.
\end{equation}
\end{theorem}
As an immediate corollary of \Cref{thm:generic-hardness} and \Cref{thm:optimal}, given the score, class $\calQ$, and full oracle families defined above, we obtain the minimax rate
\begin{equation}\label{eq:main-lower-unconstrained}
 \inf_{\calC\in\Asym}\sup_{\PP\in \calQ}\inf_{\calC^*_{\PP}\in\calO_\alpha(\PP)}
 \EE_{\PP}\ell_T(\calC,\calC^*_{\PP})
 \asymp T^{\max\{1/2,\beta\}}.
\end{equation}
The rate in \eqref{eq:main-lower-unconstrained} reflects two distinct costs.
The first is the cost of approximating an exact oracle from past data.
For split conformal prediction, disagreement with the population score oracle equals the expected error in instantaneous training-conditional coverage.
Summing \eqref{eq:iid-tcc} therefore yields an $O(\sqrt {T})$ upper bound.
The lower-bound argument for \Cref{thm:generic-hardness} shows that this order is unavoidable over $\calQ$, even without a long-run empirical coverage constraint and after taking the infimum over all exact oracles.

The second cost comes from enforcing the constraint, even if the population prediction set is already known.
As we already see from \Cref{thm:generic-hardness}, even exactly calibrated prediction sets are forced to deviate from the oracle by the long-run empirical coverage constraint, and the resulting disagreement is at least of order $T^\beta$ when $\beta\ge1/2$.

Truncation attains the combined rate by controlling both costs.
By \Cref{upper:stability}, the disagreement between the truncated baseline and the truncated oracle is at most twice the baseline's disagreement with the oracle.
To compare with the original oracle, the triangle inequality adds the oracle's own truncation cost.
Since an exact oracle has i.i.d.\ $\bern(\alpha)$ errors, \Cref{upper:bernoulli-cost} bounds this additional cost in expectation.
Thus, for any exact oracle $\calC^*_{\PP}$ and baseline $\calC^0\in\calA$, truncation gives
\begin{equation}\label{eq:generic-upper-decomposition}
\begin{split}
 \EE_{\PP}\ell_T(\Witv\circ\calC^0,\calC^*_{\PP}) \le 2 \EE_{\PP}\ell_T(\calC^0,\calC^*_{\PP})
 +K_{\alpha,\beta,\kappa}
 \begin{cases}
  1,&0\le\beta<1/2,\\
  T^\beta,&1/2\le\beta\le1.
 \end{cases}
\end{split}
\end{equation}
The bound holds for every $T\ge1$.
For split conformal prediction and the population score oracle, the first term is $O(\sqrt T)$.
When $\beta<1/2$, the additional term is bounded uniformly in $T$, so enforcing long-run empirical coverage preserves the unconstrained statistical rate.
At $\beta=1/2$, the enforcement cost becomes $\sqrt T$, matching the cost of learning the target.
For $\beta>1/2$, the enforcement cost $T^\beta$ dominates, and the optimal disagreement rate increases with the tightness of the constraint.
This leads to the sharp transition in \eqref{eq:main-lower-unconstrained}: the optimal rate remains at $\sqrt T$ up to $\beta=1/2$, and becomes $T^\beta$ beyond that threshold.

\section{Comparing to ACI with decaying step sizes}\label{sec:aci-subopt}

We next compare truncation with ACI with decaying step sizes and a fixed score function, which adjusts its threshold after each observation.
Fix $B>0$, a measurable score $s:\calX\times\calY\to[0,B]$, a deterministic initial threshold $q_1\in[0,B]$, and positive deterministic step sizes $(\eta_t)_{t\ge1}$.
ACI constructs
\begin{equation}\label{eq:aci1}
 \calC_t^{\aci}(x)=\{y:s(x,y)\le q_t\},
 \qquad
 q_{t+1}=q_t+\eta_t(\err_t^{\calC^{\aci}}-\alpha).
\end{equation}
We study the adversarial guarantee of this update and its disagreement with the closest conditionally exact oracle defined in \eqref{eq:main-oracle-calibration}.

\paragraph{Adversarial coverage guarantee.}
For nonincreasing step sizes, the coverage bound of \citet[Theorem~1]{angelopoulos2024online} gives $|\miscov_T^{\calC^{\aci}}|\le(B+\eta_1)/\eta_T$ on every sequence of scores in $[0,B]$.
In particular, for $\eta_t=ct^{-(1-\beta)}$ with $c>0$ and $\beta\in(0,1]$,
\begin{equation}\label{eq:aci-worst-case}
 \frac{|\miscov_T^{\calC^{\aci}}|}{T}
 \le\frac{B+c}{c}\,T^{-\beta}
 \qquad\text{for every }T\ge1.
\end{equation}
Thus ACI has a $(\beta,\kappa)$ adversarial guarantee with $\kappa=(B+c)/c$, indicating that, for a prescribed $\kappa>1$, choosing $c=B/(\kappa-1)$ yields the long-run coverage constraint.

\paragraph{Disagreement with the oracle.}
Now suppose the observations are i.i.d.\ under $\PP$ and the score has a density satisfying
\begin{equation}\label{eq:aci-density}
 0<\lambda_-\le f_s(u)\le\lambda_+<\infty
 \qquad\text{for almost every }u\in[0,B].
\end{equation}
Then we have the following result, where the infimum over $\calC^*_{\PP}$ below ranges over all oracles satisfying \eqref{eq:main-oracle-calibration}.
The smallest expected disagreement is achieved by the population prediction set that uses the same score, as defined in \eqref{eq:population-oracle}.

\begin{proposition}[Disagreement of ACI]\label{prop:aci-subopt}
With $\PP \in \calQ$ satisfying \eqref{eq:aci-density}, if $\eta_t\to0$ and $|\eta_{t+1}^{-1}-\eta_t^{-1}|\to0$, then
\begin{equation}\label{eq:aci-sub}
 \inf_{\calC^*_{\PP}\in\calO_\alpha(\PP)}\EE_{\PP}\bigl[\ell_T(\calC^{\aci},\calC^*_{\PP})\bigr]
 \asymp\sum_{t=1}^T\sqrt{\eta_t}.
\end{equation}
\end{proposition}

ACI's threshold fluctuates around the population quantile on the scale $\sqrt{\eta_t}$, which creates disagreement even as the step sizes decay.
Specifically, for $\eta_t=ct^{-(1-\beta)}$ with $0<\beta<1$, we have the following corollary.
\begin{corollary}\label{cor:aci-decaying-risk}
With $\PP \in \calQ$ satisfying \eqref{eq:aci-density}, for $\eta_t=ct^{-(1-\beta)}$ with $0<\beta<1$,
\begin{equation}\label{eq:aci-decaying-risk} \inf_{\calC^*_{\PP}\in\calO_\alpha(\PP)}\EE_{\PP}\bigl[\ell_T(\calC^{\aci},\calC^*_{\PP})\bigr]
 \asymp T^{(1+\beta)/2}.
\end{equation}
\end{corollary}
Since $(1+\beta)/2>\max\{1/2,\beta\}$, this grows faster than the disagreement of truncated split conformal prediction in \Cref{thm:optimal} under the same adversarial guarantee with $\kappa>1$.
Thus, with the same fixed score function, for these step sizes, ACI's cumulative disagreement exceeds the minimax rate for every $\beta \in (0,1)$.
Truncated split conformal prediction attains the optimal rate by changing the baseline set only when the constraint forces coverage or miscoverage.

\section{Numerical experiments}
\label{sec:simulations}

We study the cost of truncation and compare truncated split conformal prediction with ACI in a heteroscedastic regression model.
Independently over time, let
\[
 X_t\sim\mathrm{Unif}[0,1],\qquad
 Y_t=f(X_t)+\sigma(X_t)\vareps_t,\qquad
 \vareps_t\sim\mathrm{Unif}[-1,1],
\]
where $X_t$ and $\vareps_t$ are independent, $f(x)=\sin(2\pi x)$, and $\sigma(x)=0.2+0.3x$.
We use the fixed score $s(x,y)=\min\{|y-f(x)|/\sigma(x),1\}$, which is uniform on $[0,1]$ under this model.
We set $\alpha=0.1$ and compare all algorithms with the population prediction set $\calC^*_{\PP,t}(x)=\{y:s(x,y)\le0.9\}$.

We use the split-conformal baseline $\calC^{\cp}$ defined in \eqref{eq:conformal-baseline-set}.
Each experiment in this section uses $500$ independent replications over $T=250{,}000$ time steps.
We evaluate disagreement with the oracle using the known score distribution, which is not used to choose the conformal thresholds or update ACI.
Since the prediction sets are nested score thresholds, their disagreement probability, conditional on past observations, equals the absolute difference between their instantaneous training-conditional miscoverage probabilities.
We sum these probabilities within each replication and average the resulting totals to estimate expected cumulative disagreement.

\begin{figure}[ht]
\centering
\includegraphics[height=0.35\textwidth]{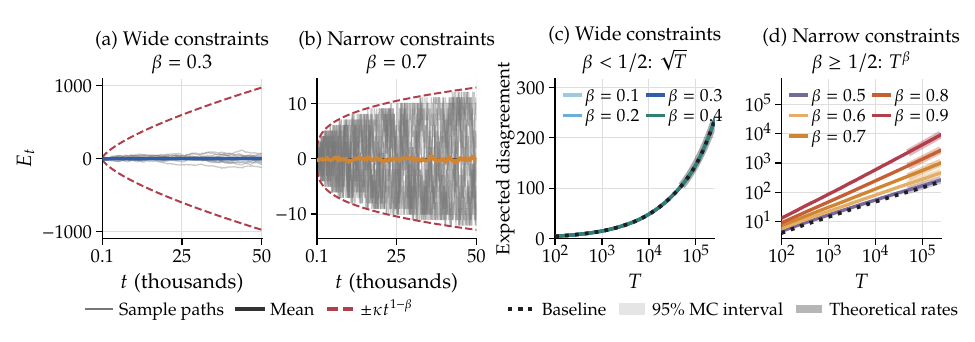}
\caption{\textbf{Long-run empirical coverage and oracle disagreement of truncated split conformal prediction.}
We use $\alpha=0.1$ and $\kappa=0.5$.
Panels~(a)--(b) show the cumulative deviations for $\beta=0.3$ and $\beta=0.7$, respectively: the first twelve paths, the mean over $500$ replications, and the bounds $\pm\kappa t^{1-\beta}$.
Paths are displayed every $100$ rounds through $t=50{,}000$.
Panels~(c)--(d) compare estimated expected cumulative disagreement with the population oracle for the truncated and unconstrained conformal algorithms.
Both panels have logarithmic horizontal axes; panel~(c) has a linear vertical axis and panel~(d) a logarithmic one.
The theoretical rates are proportional to $\sqrt T$ for $\beta<1/2$ and $T^\beta$ for $\beta\ge1/2$.}
\label{fig:sim-conformal}
\end{figure}

\paragraph{Long-run empirical coverage and oracle disagreement.}
We first truncate the split conformal baseline with $\kappa=0.5$.
The baseline continues to use all past scores, including observations from rounds in which truncation forces coverage or miscoverage.
\Cref{fig:sim-conformal} shows paths of the cumulative deviation and estimates of expected disagreement with the population oracle.
All simulated paths satisfy the long-run empirical coverage constraint at every round, up to numerical precision.
For $\beta<1/2$, disagreement remains close to that of the unconstrained baseline and follows the $\sqrt{T}$ rate predicted by theory.
For $\beta\ge1/2$, the curves are consistent with the $T^\beta$ rate in \Cref{thm:optimal}.

\paragraph{The cost of truncation.}
To measure how much truncation changes the baseline's predictions, we compare $\Witv\circ\calC^{\cp}$ with $\calC^{\cp}$ on the same observations.
Let $p_t^{\cp}=\PP(Y_t\notin\calC_t^{\cp}(X_t)\mid\calG_{t-1})$ be the baseline's instantaneous training-conditional miscoverage probability.
The training-conditional disagreement probability is $p_t^{\cp}$ when truncation forces coverage, $1-p_t^{\cp}$ when it forces miscoverage, and zero otherwise.
Summing these probabilities and averaging across replications estimates $\EE_{\PP}[\ell_T(\Witv\circ\calC^{\cp},\calC^{\cp})]$.

\begin{figure}[ht]
\centering
\includegraphics[height=0.35\textwidth]{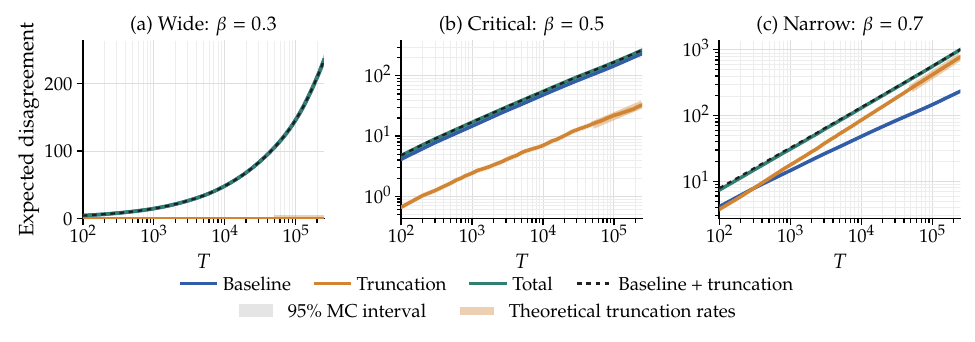}
\caption{\textbf{Baseline disagreement, truncation cost, and total disagreement.}
For $\kappa=0.5$, the three panels show $\beta=0.3,0.5,0.7$, respectively.
Each compares $\EE_{\PP}[\ell_T(\calC^{\cp},\calC^*_{\PP})]$, $\EE_{\PP}[\ell_T(\Witv\circ\calC^{\cp},\calC^{\cp})]$, and $\EE_{\PP}[\ell_T(\Witv\circ\calC^{\cp},\calC^*_{\PP})]$.
The sum curve adds the first two estimates from paired replications and bounds the corresponding total estimate from above.
All horizontal axes are logarithmic.
Panel~(a) uses a linear vertical axis to show values near zero; panels~(b)--(c) use logarithmic vertical axes.}
\label{fig:sim-decomposition}
\end{figure}

\Cref{fig:sim-decomposition} compares this truncation cost with the baseline and total disagreement for $\beta=0.3,0.5,0.7$.
The baseline disagreement is consistent with $\sqrt T$ growth in all three cases.
For $\beta=0.3$, the truncation cost remains nearly zero over the plotted horizons, so the total stays close to the baseline.
For $\beta=0.5$ and $\beta=0.7$, the truncation cost grows approximately along the $\sqrt T$ and $T^{0.7}$ rates, respectively.


\paragraph{Comparison with ACI.}
We compare ACI with truncated split conformal prediction under the same $(\beta,2)$ long-run empirical coverage guarantee, for $\beta\in\{0.3,0.5,0.7\}$.
We initialize ACI at the population threshold $q_1=0.9$ and use step sizes $\eta_t=t^{-(1-\beta)}$.
Since the score lies in $[0,1]$, \eqref{eq:aci-worst-case} gives $\kappa=2$ for ACI; we use the same value for truncation.
Within each replication, both algorithms receive the same observations and are evaluated against the same population oracle.

\begin{figure}[ht]
\centering
\includegraphics[height=0.35\textwidth]{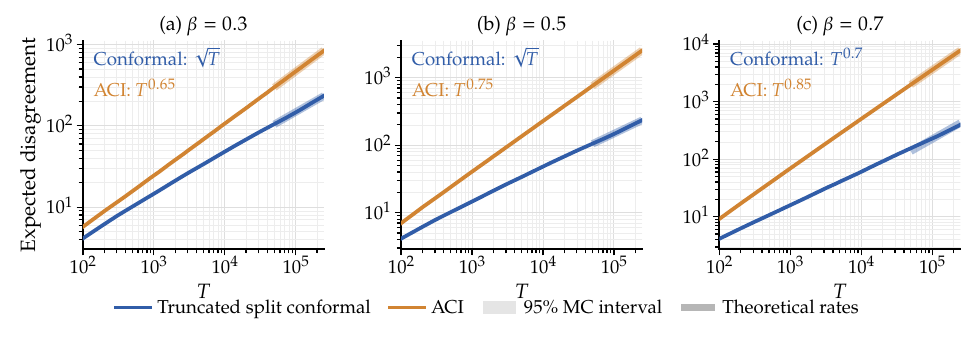}
\caption{\textbf{ACI and truncated split conformal prediction under the same long-run empirical coverage guarantee.}
Both algorithms use $\alpha=0.1$ and satisfy the $(\beta,2)$ constraint.
Each panel uses logarithmic axes and shows estimated expected cumulative disagreement with the population oracle over $500$ paired replications.
ACI uses $q_1=0.9$ and $\eta_t=t^{-(1-\beta)}$.
The theoretical rates are proportional to $T^{(1+\beta)/2}$ for ACI and $T^{\max\{1/2,\beta\}}$ for truncated split conformal prediction.}
\label{fig:sim-aci-comparison}
\end{figure}

Under this model, ACI accumulates disagreement faster than truncated split conformal prediction in all three panels of \Cref{fig:sim-aci-comparison}.
The ACI curves are consistent with the $T^{(1+\beta)/2}$ rate in \Cref{cor:aci-decaying-risk}.

\section{Related work}
\label{sec:literature}

\paragraph{Conformal prediction.}
Conformal prediction constructs prediction sets with finite-sample marginal coverage under exchangeability \citep{vovk2005algorithmic,shafer2008tutorial,angelopoulos2024theoretical}.
Inductive or split conformal prediction separates model fitting from calibration \citep{papadopoulos2002inductive,lei2018distribution}.
Conformalized quantile regression uses fitted conditional quantiles to adapt interval widths to heteroscedasticity \citep{romano2019conformalized}, while distributional conformal prediction builds scores from an estimated conditional distribution \citep{chernozhukov2021distributional}.

Extensions address covariate shift \citep{tibshirani2019conformal}, label shift \citep{podkopaev2021distribution}, and coverage loss when exchangeability fails \citep{barber2023conformal}.
\citet{cauchois2024robust} construct prediction sets with finite-sample coverage bounds for test distributions within an $f$-divergence ball.
Related work on robust risk evaluation incorporates isotonic constraints on density ratios \citep{gui2024distributionally}.
Applications include inference for counterfactuals and individual treatment effects \citep{lei2021conformal,jin2023sensitivity}, prediction with censored survival outcomes \citep{candes2023conformalized,gui2024conformalized}, missing covariates \citep{zaffran2023conformal}, and unobserved matrix entries \citep{gui2023conformalized}.

In the batch setting, \citet{areces2024two} establish a minimax lower bound of order $n^{-1/2}$ on the expected training-conditional coverage error after $n$ i.i.d.\ observations, for prediction sets defined by thresholds or intervals on a fixed score.
\citet{duchi2025few} develops sample-conditional coverage bounds for split conformal prediction and relates them to these lower bounds.
We study cumulative error in online prediction and the additional cost of a long-run empirical coverage constraint.

\paragraph{Online conformal prediction and adaptive calibration.}
Adaptive calibration uses past coverage errors to adjust prediction sets as the data distribution changes.
Adaptive conformal inference (ACI) adjusts the nominal coverage level and controls long-run empirical coverage on arbitrary data sequences \citep{gibbs2021adaptive}.
Subsequent methods adapt learning rates with local dynamic-regret guarantees \citep{gibbs2024conformal}, or aggregate ACI predictors with different learning rates \citep{zaffran2022adaptive}.
Betting-based methods use parameter-free online optimization to avoid choosing a learning rate \citep{podkopaev2024adaptive}.

Direct threshold updates connect long-run empirical coverage control with quantile estimation.
Conformal PID control combines quantile tracking, correction for errors in cumulative coverage, and forecasts of future score quantiles \citep{angelopoulos2023conformal}.
For fixed bounded scores, appropriately decaying step sizes preserve long-run empirical coverage and yield quantile convergence under i.i.d.\ sampling and a unique target quantile \citep{angelopoulos2024online}.

Time-series methods use past residuals and forecasts to adapt prediction intervals to temporal dependence.
EnbPI uses bootstrap ensemble residuals and updates their empirical quantiles as new observations arrive \citep{xu2021conformal}.
Sequential predictive conformal inference instead estimates conditional quantiles of future residuals from past residuals \citep{xu2023sequential}.
Bellman conformal inference uses multi-step forecasts and dynamic programming to choose prediction intervals \citep{yang2024bellman}.
\citet{jiang2026leave} document severe undercoverage of the leave-one-out jackknife under short-range temporal dependence.
Their leave-a-window-out modification omits blocks of observations and gives approximate marginal-coverage guarantees under stability of the fitting algorithm and quantitative bounds on temporal dependence.

\cite{srinivas2026online} studies the trade-off between empirical coverage and average interval length, measured relative to the shortest fixed interval that achieves the target coverage in hindsight.
\citet{areces2025online} formulate calibration as an online optimization problem, obtaining the long-run empirical coverage in the adversarial setting and convergence of instantaneous training-conditional coverage under assumptions that include a linear model for conditional score quantiles.
When the data law changes, a central goal is to keep training-conditional coverage close to its target.
\citet{liang2026optimal} study $\inst$, defined in \eqref{eq:instant}, for independent observations with unknown change points or bounded cumulative distribution drift.
Their methods combine drift detection and recalibration and attain minimax rates up to logarithmic factors within the algorithm classes they study, but do not impose a long-run empirical coverage constraint that holds on every data sequence.

\paragraph{Regret, local coverage, and group coverage.}
Strongly adaptive online methods bound quantile-loss regret over all time intervals \citep{bhatnagar2023improved}, where, however, quantile-loss regret alone does not directly imply coverage guarantees.
Multivalid methods control empirical coverage across groups and threshold ranges under smoothness conditions on the score distributions \citep{bastani2022practical}.
\citet{ramalingam2025relationship} distinguish regret from coverage, relate swap regret to multivalidity, and develop group-conditional extensions of ACI.
Other work studies exponentially discounted regret and coverage \citep{zhang2024discounted}, selection among candidate prediction models \citep{hajihashemi2024multi}, and gradient equilibrium, which relates vanishing average loss gradients to online calibration \citep{angelopoulos2025gradient}.


\section{Discussion}
\label{sec:discussion}

We seek a \emph{best-of-both-worlds} online algorithm whose prediction sets satisfy the long-run empirical coverage constraint \eqref{eq:corridor_intro} on every data sequence and disagree as little as possible with oracles whose instantaneous training-conditional coverage is exactly $1-\alpha$ for i.i.d. data.
We show that, over the class $\calQ$ of i.i.d.\ distributions with a nonatomic score, the minimax expected disagreement with the closest such oracle is of order $T^{\max\{1/2,\beta\}}$ under this constraint.
The constraint therefore leaves the unconstrained rate $\sqrt{T}$ unchanged for $\beta\le1/2$ and raises it to $T^\beta$ for $\beta>1/2$.
Truncated split conformal prediction attains this rate.
With the same fixed score, ACI with step sizes proportional to $t^{-(1-\beta)}$ satisfies the same constraint, but its expected disagreement is of order $T^{(1+\beta)/2}$ when the score has a density bounded away from zero and infinity.
This rate exceeds $T^{\max\{1/2,\beta\}}$ for every $0<\beta<1$.
This paper also leaves the following questions open.
Since the minimax-optimal construction may not be unique, it may be of interest to refine the optimality criterion by quantifying the accuracy of such prediction sets via the volume. 
Open questions include how to retain informative prediction sets and adaptive instantaneous training-conditional coverage under distribution drift while preserving a prescribed adversarial constraint.

\section*{Acknowledgement}
R.F.B. was partially supported by the Office of Naval Research via grant N00014-24-1-2544.
C.M. was partially supported by the National Science Foundation via the CAREER Award DMS2443867.

\section*{Disclosure of AI usage}
GPT-6-Astra and Claude-Fable-5.1 models were used to produce numerical experiments, polish writing, suggest literature, and audit existing proofs; the authors take full responsibility for the contents of the paper.

\bibliography{main}
\bibliographystyle{apalike}

\newpage
\appendix

\addtocontents{toc}{\protect\setcounter{tocdepth}{3}}
\tableofcontents

\section{Results under asymmetric constraints}
\label{app:asymmetric}

The main text uses the symmetric constraint \eqref{eq:corridor_intro}.
Here we allow its two sides to grow at different rates and prove the main results in this setting.
The symmetric statements follow by taking $b=d=\beta$ and $c_-=c_+=\kappa$.
The ACI result is proved separately in Appendix~\ref{app:aci}.

\subsection{Setting and notation}
\label{app:asymmetric-setting}

Fix $\alpha\in(0,1)$, $b,d\in[0,1]$, and $c_-,c_+>0$ with $c_-+c_+\ge1$.
Write
\[
 L_t=c_-t^{1-d},\qquad U_t=c_+t^{1-b},\qquad
 \bunder=\min\{b,d\},\qquad \bbar=\max\{b,d\},
 \qquad t\ge1.
\]
Set $L_0=U_0=0$.

\paragraph{Notations.}
We write $\sigma=\sqrt{\alpha(1-\alpha)}$, $\nu=\min\{\alpha,1-\alpha\}$, and $x_+=\max\{x,0\}$.
Write $Z_t=(X_t,Y_t)\in\calZ=\calX\times\calY$ and $\calG_{t-1}=\sigma((Z_r)_{r<t})$.
For a fixed measurable score $s$, put $s_t=s(Z_t)$.
Write $\tilde{\calQ}$ for all i.i.d.\ distributions on the observation sequence and $\calQ$ for the nonempty subclass under which $s_1$ is atomless.
For an i.i.d.\ distribution $\PP$, denote its one-observation marginal by $\mu_\PP$, so $\PP=\mu_\PP^{\otimes\mathbb N}$.
For fixed $\PP$ and $s$, abbreviate the score CDF $F_{\PP,s}$ to $F$ and retain $q_{\PP,s}$ for the population quantile in \eqref{eq:population-oracle}.

Write $e^{\calC}=(\err_t^{\calC})_{t\ge1}$ for the error sequence of $\calC$ and $p_t^{\calC}=\PP(\err_t^{\calC}=1\mid\calG_{t-1})$ for its instantaneous training-conditional miscoverage probability.
The error on a history $h$ and fresh observation $z$ is a product-measurable function $e_t(h,z)$.
Use $e^0$ and $e^*$ for the error sequences of the baseline and comparison oracle, respectively.
When a single algorithm is fixed, we may omit its superscript.
For binary paths $u,v$, write
\[
 \miscov_t(u)=\sum_{r=1}^t(u_r-\alpha),\qquad
 D_T(u,v)=\sum_{r=1}^T|u_r-v_r|,\qquad \miscov_0(u)=0 .
\]
Thus $\miscov_t^{\calC}=\miscov_t(e^{\calC})$ and $\ell_T(\calC,\calC')=D_T(e^{\calC},e^{\calC'})$.
When comparing a constrained path $u$ to a reference path $v$, abbreviate their centered sums to $\miscov_t=\miscov_t(u)$ and $S_t=\miscov_t(v)$, respectively.

\paragraph{The coverage constraint.}
The class $\Aasym$ consists of algorithms whose centered errors satisfy
\begin{equation}\label{eq:asymmetric-constraint}
 -L_t\le\miscov_t^{\calC}\le U_t
 \qquad\text{for every }t\ge1\text{ and every data sequence}.
\end{equation}
Thus $b$ controls the upper side and $d$ the lower side.
All algorithms and oracles in this appendix are deterministic, as in the main text.

Every binary path $u$ satisfies $-\alpha t\le\miscov_t(u)\le(1-\alpha)t$.
The first horizon at which the constraint excludes a binary sequence is
\begin{equation}\label{eq:asym-first-active}
 T_0=\inf\{t\in\mathbb N_{\ge1}: U_t<(1-\alpha)t\ \text{or}\ L_t<\alpha t\},
 \qquad\inf\varnothing=\infty.
\end{equation}
We call the constraint nontrivial when $T_0<\infty$.
Otherwise $b=d=0$, $c_+\ge1-\alpha$, and $c_-\ge\alpha$, so $\Aasym=\calA$.

For nontrivial constraints, the cost of imposing the constraint on an exact oracle has three regimes:
\begin{equation}\label{eq:asym-cost}
 g_{b,d}(T)=
 \begin{cases}
  1,&\bbar<1/2,\\
  \sqrt T,&\bunder\le1/2\le\bbar,\\
  T^{\bunder},&\bunder>1/2.
 \end{cases}
\end{equation}
Both boundaries grow faster than $\sqrt t$ in the first regime.
At least one boundary grows at most as fast as $\sqrt t$ in the second, and both grow more slowly in the third.
In the symmetric case, $g_{\beta,\beta}(T)$ is $1$ for $\beta<1/2$ and $T^\beta$ for $\beta\ge1/2$.

\paragraph{Asymmetric truncation.}
For a baseline $\calC^0$, let $\widehat\calC=\Gamma_{\rm asym}\circ\calC^0$ denote the truncated algorithm for the fixed parameters $(b,d,c_-,c_+)$, which is defined by
\begin{equation}\label{eq:asym-truncation}
 \widehat\calC_t(x)=
 \begin{cases}
  \calY,&\miscov_{t-1}^{\widehat\calC}+1-\alpha>U_t,\\
  \varnothing,&\miscov_{t-1}^{\widehat\calC}-\alpha<-L_t,\\
  \calC_t^0(x),&\text{otherwise},
 \end{cases}
 \qquad\miscov_0^{\widehat\calC}=0.
\end{equation}
For a binary input $e$, write $\hat e=\Gamma_{\rm asym}(e)$ for the output on errors.
The first two conditions in \eqref{eq:asym-truncation}, evaluated at the centered sum of $\hat e$, specify rounds forced from above and below, respectively.
In particular,
\begin{equation}\label{upper:error-map}
 e^{\Gamma_{\rm asym}\circ\calC^0}=\Gamma_{\rm asym}(e^{\calC^0}),\qquad
 \hat e_t=
 \begin{cases}
  0,&\text{forced from above},\\
  1,&\text{forced from below},\\
  e_t,&\text{otherwise}.
 \end{cases}
\end{equation}
The two conditions are disjoint since $L_t+U_t\ge1$.
The baseline uses every observation, including observations from forced rounds.
Under $b=d=\beta$ and $c_-=c_+=\kappa$, $\Gamma_{\rm asym}$ specializes to $\Witv$ in \eqref{eq:intervention}.

Let $N_T^+,N_T^-$ count the rounds forced from above and below through $T$, and $M_T^+,M_T^-$ the rounds on which those respective branches change the input error.
Write $N_T=N_T^++N_T^-$ and use subscript $\infty$ for the corresponding limits.
Appendix~\ref{upper:section} proves feasibility and stability.

\paragraph{Oracle risks.}
When $\PP$ is fixed, write $\calC^*_{\PP}$ for an oracle in the exact-oracle family $\calO_\alpha(\PP)$ of \eqref{eq:main-oracle-calibration}.
Define
\begin{equation}\label{eq:asym-risk-notation}
 \begin{aligned}
 r_T(\calC;\PP)&=\inf_{\calC^*_{\PP}\in\calO_\alpha(\PP)}\EE_\PP\ell_T(\calC,\calC^*_{\PP}),\\
 R_T(\calC)&=\sup_{\PP\in\calQ}r_T(\calC;\PP),\qquad
 R_T^0=\inf_{\calC\in\calA}R_T(\calC).
 \end{aligned}
\end{equation}
All constants in the constraint bounds below depend only on $\alpha,b,d,c_-,c_+$.

\subsection{Unconstrained minimax rate}
\label{app:unconstrained}

Without a long-run empirical coverage constraint, the cost of learning from past observations is of order $\sqrt T$, even when the comparison ranges over all exact oracles.

\begin{lemma}\label{lem:learning-full}
For every $T\ge1$,
\[
 \frac{\nu^2}{32}\sqrt T\le R_T^0\le3\sqrt T.
\]
\end{lemma}

\noindent Appendix~\ref{learning:section} proves the lemma and shows that split conformal prediction attains its upper bound.

\subsection{Lower bounds}

\begin{theorem}[Disagreement forced by an asymmetric constraint]\label{thm:asym-lower}
Suppose that $T_0<\infty$.
For every i.i.d.\ law $\PP$ admitting an exact oracle, every such oracle $\calC^*_{\PP}$, and every online algorithm $\calC$ satisfying \eqref{eq:asymmetric-constraint},
\begin{equation}\label{eq:asym-pointwise-lower}
 \EE_\PP\ell_T(\calC,\calC^*_{\PP})\ge c\,g_{b,d}(T)
 \qquad(T\ge T_1),
\end{equation}
where $c>0$ and $T_1$ are uniform over the law, algorithm, and oracle.
This statement allows $\calC$ to know $\PP$.
Consequently, every $\calC\in\Aasym$ satisfies
\begin{equation}\label{eq:asym-lower-decomposition}
 R_T(\calC)\gtrsim R_T^0+g_{b,d}(T).
\end{equation}
\end{theorem}
Appendix~\ref{app:lower} proves the constraint bound and combines the two costs.
Under the symmetric specialization, this gives \Cref{thm:generic-hardness}.

\subsection{Truncation and the minimax rate}

\begin{theorem}[Stability of truncation]\label{thm:asym-upper}
For every baseline $\calC^0\in\calA$, the algorithm $\Gamma_{\rm asym}\circ\calC^0$ belongs to $\Aasym$.
For every i.i.d.\ law $\PP$ and every exact oracle $\calC^*_{\PP}$ under that law,
\begin{equation}\label{eq:asym-transfer}
 \EE_\PP\ell_T(\Gamma_{\rm asym}\circ\calC^0,\calC^*_{\PP})
 \le2\EE_\PP\ell_T(\calC^0,\calC^*_{\PP})+K\,g_{b,d}(T)
 \qquad(T\ge1),
\end{equation}
where $K<\infty$ is independent of the baseline, law, and oracle.
\end{theorem}

This result does not require a continuous score distribution or calibrated baseline errors.
It uses only the exact oracle condition of the comparison algorithm.
Appendix~\ref{upper:section} proves it, and the symmetric case gives \eqref{eq:generic-upper-decomposition}.

\begin{theorem}[Minimax disagreement under asymmetric constraints]\label{thm:asym-optimal}
For the fixed score $s$, nonempty class $\calQ$, and full exact-oracle families of the main text,
\begin{equation}\label{eq:asym-optimal}
 \inf_{\calC\in\Aasym}R_T(\calC)
 \asymp T^{\max\{1/2,\bunder\}}.
\end{equation}
The single algorithm $\Gamma_{\rm asym}\circ\calC^{\cp}$ attains the upper bound.
\end{theorem}

Appendix~\ref{app:optimal} combines the learning and constraint bounds to prove this theorem.
Taking $b=d=\beta$ and $c_-=c_+=\kappa$ gives \Cref{thm:optimal}.

\subsection{Proof of Theorem~\ref{thm:asym-optimal}}
\label{app:optimal}

\begin{proof}
We first compare the three possible constraint costs with the learning cost.
We then establish matching lower and upper bounds.

\paragraph{Step 1: identify the resulting rate.}
The definition of $g_{b,d}$ gives
\begin{equation}\label{optimal:rates}
 \sqrt T+g_{b,d}(T)\asymp T^{\max\{1/2,\bunder\}}.
\end{equation}
Indeed, $g_{b,d}(T)$ is bounded when $\bbar<1/2$ and is $\sqrt T$ when $\bunder\le1/2\le\bbar$.
In both cases the sum has order $\sqrt T$.
When $\bunder>1/2$, the term $T^{\bunder}$ is the larger term.

\paragraph{Step 2: bound every constrained algorithm from below.}
Suppose first that the constraint is nontrivial.
By \Cref{thm:asym-lower} and \Cref{lem:learning-full}, every $\calC\in\Aasym$ satisfies
\[
 R_T(\calC)\gtrsim\sqrt T+g_{b,d}(T)
 \asymp T^{\max\{1/2,\bunder\}}.
\]
The constants and initial horizon are uniform over $\calC$, so taking the infimum over $\Aasym$ preserves this bound.
If the constraint is trivial, then $\Aasym=\calA$ and $\bunder=0$.
In that case \Cref{lem:learning-full} gives the same conclusion directly.

\paragraph{Step 3: attain the rate with truncated conformal prediction.}
Take $\calC^0=\calC^{\cp}$ in \eqref{eq:asym-truncation}.
By \Cref{thm:asym-upper}, this algorithm belongs to $\Aasym$ on every data sequence.
For each $\PP\in\calQ$, the population score oracle is exact, so
\begin{align*}
 \inf_{\calC^*_{\PP}\in\calO_\alpha(\PP)}
   \EE_\PP\ell_T(\widehat\calC,\calC^*_{\PP})
 &\le\EE_\PP\ell_T(\widehat\calC,\calC^*_{\PP,s})\\
 &\le2\EE_\PP\ell_T(\calC^{\cp},\calC^*_{\PP,s})+K g_{b,d}(T)\\
 &\le6\sqrt T+K g_{b,d}(T).
\end{align*}
The last bound is uniform over $\PP$, and \eqref{optimal:rates} gives the desired upper rate.
Neither the conformal baseline nor truncation depends on $T$, so the same algorithm works for all horizons.

Finally, setting $b=d=\beta$ and $c_-=c_+=\kappa$ turns $\Gamma_{\rm asym}$ into the truncation rule \eqref{eq:intervention}.
The upper and lower bounds then give \eqref{eq:main-lower-unconstrained} and \eqref{eq:rates}.
\end{proof}




\section{Proof of the asymmetric lower bound}
\label{app:lower}

The proof is based on the following facts.

\begin{lemma}\label{lower:oracle-iid}
If a deterministic online oracle satisfies \eqref{eq:main-oracle-calibration}, its errors are i.i.d.\ $\bern(\alpha)$.
\end{lemma}

\begin{lemma}[Comparing a constrained path with a reference path]\label{lower:path-comparison}
Suppose that $-L_t\le \miscov_t\le U_t$ for every $t\le T$.
For $0\le r<r'\le T$,
\begin{equation}\label{lower:block}
 \sum_{t=r+1}^{r'}|e_t-e_t^*|
 \ge |S_{r'}-S_r|-(c_-+c_+)T^{1-\bunder}.
\end{equation}
At the terminal horizon,
\begin{equation}\label{lower:terminal}
 D_T(e,e^*)\ge(S_T-U_T)_+,\qquad D_T(e,e^*)\ge(-S_T-L_T)_+.
\end{equation}
\end{lemma}

\begin{lemma}[Fluctuations of the Bernoulli walk]\label{lower:walk}
Suppose the reference errors are i.i.d.\ $\bern(\alpha)$.
There is $a_\alpha^{\rm rw}>0$ such that $\EE|S_n|\ge a_\alpha^{\rm rw}\sqrt n$ for every $n\ge1$.
For every fixed $a>0$,
\[
 \lim_{n\to\infty}\PP(S_n\ge a\sqrt n)
 =\lim_{n\to\infty}\PP(-S_n\ge a\sqrt n)
 =1-\Phi(a/\sigma)>0,
\]
where $\Phi$ is the standard normal distribution function.
\end{lemma}

The first two lemmas are elementary facts specific to our comparison.
The last uses the classical central limit theorem and convergence of moments under uniform integrability; we recall the short application below rather than prove these standard results.

\subsection{Proof of Theorem~\ref{thm:asym-lower}}

\begin{proof}
Fix $\PP$, $\calC$, and $\calC^*_{\PP}$ as in the theorem.
Apply \Cref{lower:path-comparison} to $e=e^{\calC}$ and $e^*=e^{\calC^*_{\PP}}$.
By \Cref{lower:oracle-iid}, $S_t$ has the law in \Cref{lower:walk}, independently of the choices of law, algorithm, and oracle.

\paragraph{Step 1: two tight boundaries.}
Suppose $\bunder>1/2$ and write $C=c_-+c_+$.
On any block, the constrained path can move by at most $CT^{1-\bunder}$.
Choose blocks long enough that the oracle's expected displacement is at least twice this amount:
\[
 A=\left(\frac{2C}{a_\alpha^{\rm rw}}\right)^2,\qquad
 n=\left\lceil AT^{2-2\bunder}\right\rceil,\qquad m=\lfloor T/n\rfloor.
\]
For each of the $m$ complete blocks, \Cref{lower:path-comparison,lower:walk} give expected disagreement at least $a_\alpha^{\rm rw}\sqrt n/2$.
Thus
\[
 \EE D_T(e,e^*)\ge \frac{a_\alpha^{\rm rw}}{2}m\sqrt n.
\]
Since $2-2\bunder<1$, we have $n\le T/2$ for all sufficiently large $T$, including when $\bunder=1$.
Then $m\ge T/(2n)$ and $n\le(A+1)T^{2-2\bunder}$, giving
\begin{equation}\label{lower:narrow-rate}
 \EE D_T(e,e^*)\ge\frac{a_\alpha^{\rm rw} T}{4\sqrt n}
 \ge\frac{a_\alpha^{\rm rw}}{4\sqrt{A+1}}T^{\bunder}.
\end{equation}

\paragraph{Step 2: at least one tight boundary.}
Suppose $\bunder\le1/2\le\bbar$.
If $b\ge1/2$, then $U_T\le c_+\sqrt T$.
On $\{S_T\ge(c_++1)\sqrt T\}$, \eqref{lower:terminal} gives $D_T(e,e^*)\ge\sqrt T$.
Consequently,
\[
 \EE D_T(e,e^*)\ge\sqrt T\,\PP\{S_T\ge(c_++1)\sqrt T\}\gtrsim\sqrt T
\]
by \Cref{lower:walk}.
If $d\ge1/2$, apply the same argument to $-S_T$ and $c_-$.

\paragraph{Step 3: two wide boundaries.}
Suppose $\bbar<1/2$.
Since $T_0<\infty$, at time $T_0$ the constraint excludes either the all-one or all-zero error path.
The oracle follows these paths with probabilities $\alpha^{T_0}$ and $(1-\alpha)^{T_0}$, respectively.
On the excluded path, the constrained algorithm must disagree at least once.
As $D_T(e,e^*)$ is nondecreasing,
\[
 \EE D_T(e,e^*)\ge\nu^{T_0}>0\qquad(T\ge T_0).
\]
Together, Steps~1--3 prove \eqref{eq:asym-pointwise-lower}, with constants and starting horizons depending only on the constraint parameters and $\alpha$.

\paragraph{Step 4: add the learning cost.}
Now let $\calC\in\Aasym$ be distribution-independent.
The pointwise bound survives the infimum over exact oracles and the supremum over $\calQ$, so $R_T(\calC)\ge c g_{b,d}(T)$ for large $T$.
Also $R_T(\calC)\ge R_T^0$ because $\calC\in\calA$.
Taking the larger of these two bounds gives
\[
 R_T(\calC)\ge\frac{\min\{1,c\}}{2}\bigl(R_T^0+g_{b,d}(T)\bigr).
\]
Finally, \Cref{lem:learning-full} gives $R_T^0\asymp\sqrt T$ for the class and oracle family of the main text.
This proves the theorem.
\end{proof}

\subsection{Proof of Lemma~\ref{lower:oracle-iid}}

\begin{proof}
Because the oracle is deterministic and online, $e^*_1,\ldots,e^*_{t-1}$ are measurable with respect to $\calG_{t-1}$.
The exact oracle condition gives $\PP(e^*_t=1\mid\calG_{t-1})=\alpha$ almost surely.
For any $a_1,\ldots,a_t\in\{0,1\}$, conditioning on $\calG_{t-1}$ therefore gives
\[
 \PP(e^*_1=a_1,\ldots,e^*_t=a_t)
 =\alpha^{a_t}(1-\alpha)^{1-a_t}
   \PP(e^*_1=a_1,\ldots,e^*_{t-1}=a_{t-1}).
\]
Iteration proves the product Bernoulli law.
In fact, this argument only needs exact instantaneous training-conditional coverage, not independence of the original observations.
\end{proof}

\subsection{Proof of Lemma~\ref{lower:path-comparison}}

\begin{proof}
For $0\le r<r'\le T$, the triangle inequality gives
\[
 \sum_{t=r+1}^{r'}|e_t-e_t^*|
 \ge\left|(S_{r'}-S_r)-(\miscov_{r'}-\miscov_r)\right|.
\]
Both $\miscov_r$ and $\miscov_{r'}$ lie in $[-c_-T^{1-\bunder},c_+T^{1-\bunder}]$; this includes $\miscov_0=0$.
Their difference has absolute value at most $(c_-+c_+)T^{1-\bunder}$, proving \eqref{lower:block}.
At $r=0$, the bounds $\miscov_T\le U_T$ and $\miscov_T\ge-L_T$ give $D_T(e,e^*)\ge S_T-U_T$ and $D_T(e,e^*)\ge-S_T-L_T$.
Since $D_T(e,e^*)\ge0$, the positive-part bounds follow.
\end{proof}

\subsection{Justification of Lemma~\ref{lower:walk}}

The classical i.i.d.\ central limit theorem gives $S_n/\sqrt n\Rightarrow N(0,\sigma^2)$.
Also $\EE(S_n/\sqrt n)^2=\sigma^2$, so $|S_n|/\sqrt n$ is uniformly integrable.
Convergence of first moments therefore gives
\[
 \frac{\EE|S_n|}{\sqrt n}\longrightarrow\sigma\sqrt{2/\pi}>0;
\]
see \citet[Chapter~2, Theorem~2.20]{van2000asymptotic} for these standard convergence results.
Every finite $n$ also has $\EE|S_n|>0$, so the infimum of these normalized expectations is positive.
Finally, convergence of the two tail probabilities follows from the central limit theorem because the normal distribution has no atoms.

\section{Proof of the truncation upper bound}
\label{upper:section}

This section proves \Cref{thm:asym-upper}, and we first present the following lemmas.
In the notation of Appendix~\ref{app:asymmetric-setting}, take $(u,v)=(\hat e,e)$, the truncated output and its input.

\begin{lemma}[Feasibility and stability]\label{upper:stability}
For every binary input $e$, truncation satisfies
\[
 -L_t\le \miscov_t\le U_t\qquad\text{for every }t\ge1.
\]
For every two binary inputs $e,e'$ and every $T\ge1$,
\begin{equation}\label{upper:stability-bound}
 D_T\bigl(\Gamma_{\rm asym}(e),\Gamma_{\rm asym}(e')\bigr)\le2D_T(e,e').
\end{equation}
\end{lemma}

\begin{lemma}[Expected cost on Bernoulli errors]\label{upper:bernoulli-cost}
Let the input errors $e_t$ be i.i.d.\ $\bern(\alpha)$.
There is a finite constant $K$ such that
\begin{equation}\label{upper:bernoulli-bound}
 \EE D_T\bigl(\Gamma_{\rm asym}(e),e\bigr)\le\EE N_T\le K g_{b,d}(T)
 \qquad\text{for every }T\ge1.
\end{equation}
If $\bbar<1/2$, then $\EE N_\infty<\infty$.
\end{lemma}

\Cref{upper:stability} shows that truncating two nearby error sequences keeps them nearby.
\Cref{upper:bernoulli-cost} bounds the additional disagreement needed to truncate an exact oracle, independently of its data law.
We first combine these facts and then prove the lemmas.

\subsection{Proof of Theorem~\ref{thm:asym-upper}}
\label{upper:main-proof}

\begin{proof}
Compare the baseline and oracle on the same observations.

\paragraph{Step 1: enforce the constraint.}
The first part of \Cref{upper:stability} applies to every realization of $e^0$.
Hence $\Gamma_{\rm asym}\circ\calC^0$ satisfies the asymmetric coverage constraint on every data sequence.

\paragraph{Step 2: compare the two truncated sequences.}
By the triangle inequality and \eqref{upper:stability-bound},
\begin{align*}
 \ell_T(\Gamma_{\rm asym}\circ\calC^0,\calC^*_{\PP})
 &=D_T\bigl(\Gamma_{\rm asym}(e^0),e^*\bigr)\\
 &\le D_T\bigl(\Gamma_{\rm asym}(e^0),\Gamma_{\rm asym}(e^*)\bigr)
       +D_T\bigl(\Gamma_{\rm asym}(e^*),e^*\bigr)\\
 &\le2\ell_T(\calC^0,\calC^*_{\PP})
       +D_T\bigl(\Gamma_{\rm asym}(e^*),e^*\bigr).
\end{align*}
This comparison is pathwise and requires no independence between the baseline and the oracle.

\paragraph{Step 3: bound the oracle's truncation cost.}
By \Cref{lower:oracle-iid}, $e^*$ is i.i.d.\ $\bern(\alpha)$.
Taking expectations and applying \Cref{upper:bernoulli-cost} gives
\[
 \EE_{\PP}\ell_T(\Gamma_{\rm asym}\circ\calC^0,\calC^*_{\PP})
 \le2\EE_{\PP}\ell_T(\calC^0,\calC^*_{\PP})+K g_{b,d}(T)
 \qquad(T\ge1).
\]
The constant is independent of the baseline, the oracle, and the law, as required.
\end{proof}

\subsection{Proof of Lemma~\ref{upper:stability}}
\label{upper:stability-proof}

\begin{proof}
We prove the lemma in the following steps.

\paragraph{Step 1: write the allowed error counts.}
Let $K_t=\sum_{r\le t}\hat e_r$, with $K_0=0$, and define
\[
 a_t=\lceil\alpha t-L_t\rceil,
 \qquad z_t=\lfloor\alpha t+U_t\rfloor
 \qquad(t\ge1),\qquad a_0=z_0=0.
\]
The coverage constraint is equivalent to $a_t\le K_t\le z_t$.
The interval $[\alpha t-L_t,\alpha t+U_t]$ has length $L_t+U_t\ge c_-+c_+\ge1$, so $a_t\le z_t$.
Also, for every $t\ge1$,
\begin{equation}\label{upper:count-movement}
 z_t\ge z_{t-1},\qquad a_t\le a_{t-1}+1.
\end{equation}
For $t\ge2$, these inequalities follow because $U_t$ is nondecreasing and $\alpha t-L_t$ increases by at most $\alpha<1$ in one round.
For $t=1$, they follow from $z_1\ge0$ and $a_1\le1$.

\paragraph{Step 2: identify truncation with clipping.}
Define $\operatorname{clip}_t(x)=\min\{z_t,\max\{a_t,x\}\}$.
We prove by induction that
\begin{equation}\label{upper:clipping}
 K_t=\operatorname{clip}_t(K_{t-1}+e_t),\qquad a_t\le K_t\le z_t.
\end{equation}
Assume $a_{t-1}\le K_{t-1}\le z_{t-1}$.
The upper trigger in \eqref{eq:asym-truncation} is equivalent to $K_{t-1}+1>z_t$.
By \eqref{upper:count-movement}, this can occur only when $K_{t-1}=z_t$ and forcing error zero then gives the clipped value.
The lower trigger is equivalent to $K_{t-1}<a_t$.
Again, by \eqref{upper:count-movement}, this can occur only when $K_{t-1}+1=a_t$, which forces the error to be 1, giving the clipped value.
Both triggers cannot hold simultaneously, since it would imply $z_t=K_{t-1}<a_t$.
If neither holds, both $K_{t-1}$ and $K_{t-1}+1$ lie in $[a_t,z_t]$, and the input error is unchanged.
This proves \eqref{upper:clipping} and feasibility.

\paragraph{Step 3: change one input at a time.}
Suppose first that $e$ and $e'$ differ only at round $r$, with $e_r=1$ and $e'_r=0$.
Let $K_t,K'_t$ be their clipped counts.
The function $\operatorname{clip}_t$ is nondecreasing and satisfies
\[
 0\le\operatorname{clip}_t(x)-\operatorname{clip}_t(y)\le x-y
 \qquad(x\ge y).
\]
Thus $K_t-K'_t$ is zero before $r$, is either zero or one at $r$, and is nonnegative and nonincreasing afterward because the later inputs agree.
The gap, which takes integer values, can change at most twice: from zero to one at $r$, and from one back to zero at a later round.
The difference of the output errors equals the increment of this gap, so $D_T(\Gamma_{\rm asym}(e),\Gamma_{\rm asym}(e'))\le2$.
The case $e_r=0, e'_r=1$ follows by exchanging the two inputs. 
For arbitrary inputs, replace their differing coordinates in the first $T$ rounds, one at a time, and apply the triangle inequality.
The first $T$ outputs depend only on those coordinates, giving \eqref{upper:stability-bound}.
\end{proof}

\subsection{Proof of Lemma~\ref{upper:bernoulli-cost}}
\label{upper:cost-proof}

\paragraph{Counts and three bounds.}

For the i.i.d.\ $\bern(\alpha)$ input in \Cref{upper:bernoulli-cost}, the counts of Appendix~\ref{app:asymmetric-setting} give
\begin{equation}\label{upper:corrections}
 D_T(\Gamma_{\rm asym}(e),e)=M_T^++M_T^-\le N_T,
 \qquad \miscov_t=S_t-M_t^++M_t^-.
\end{equation}
The trigger at round $t$ depends only on $e_1,\ldots,e_{t-1}$ and is independent of $e_t$.
It follows that
\begin{equation}\label{upper:expected-changes}
 \EE M_T^+=\alpha\EE N_T^+,
 \qquad \EE M_T^-=(1-\alpha)\EE N_T^-.
\end{equation}
We use three bounds, depending on whether a boundary grows faster than $\sqrt{t} $.

\begin{lemma}\label{upper:loose}
If $d<1/2$, then $\EE N_\infty^-<\infty$.
If $b<1/2$, then $\EE N_\infty^+<\infty$.
\end{lemma}

\begin{lemma}\label{upper:mixed}
If $d<1/2\le b$, then $\EE N_T^+\le C\sqrt T$ for every $T\ge1$.
If $b<1/2\le d$, then $\EE N_T^-\le C\sqrt T$ for every $T\ge1$.
\end{lemma}

\begin{lemma}\label{upper:tight}
If $\bunder\ge1/2$, then $\EE N_T\le C T^{\bunder}$ for every $T\ge1$.
\end{lemma}

\begin{proof}[Proof of \Cref{upper:bernoulli-cost}]
The first inequality in \eqref{upper:bernoulli-bound} follows from \eqref{upper:corrections}.
\begin{itemize}
    \item If $\bbar<1/2$, \Cref{upper:loose} bounds $\EE N_\infty^+$ and $\EE N_\infty^-$, so $\EE N_T$ is bounded uniformly in $T$.
    \item If $\bunder<1/2\le\bbar$, \Cref{upper:loose,upper:mixed} give $\EE N_T\le C\sqrt T$.
    \item If $\bunder=1/2$, \Cref{upper:tight} gives the same bound.
    \item Finally, if $\bunder>1/2$, \Cref{upper:tight} gives $\EE N_T\le C T^{\bunder}$.
\end{itemize}
These are exactly the three rates in $g_{b,d}(T)$.
\end{proof}

\subsection{Proof of Lemma~\ref{upper:loose}}

\begin{proof}
Assume $d<1/2$ and that round $t$ is forced from below.
Let $r<t$ be the most recent upper forced round, or set $r=0$ if there is none.
If $r>0$, its upper trigger gives $\miscov_{r-1}+1-\alpha>U_r\ge0$, and hence $\miscov_r=\miscov_{r-1}-\alpha>-1$.
If $r=0$, the same lower bound holds because $\miscov_0=0$.
There are no upper changes between $r+1$ and $t-1$, so \eqref{upper:corrections} gives
\[
 \miscov_{t-1}\ge \miscov_r+S_{t-1}-S_r>-1+S_{t-1}-S_r.
\]
The lower trigger gives $\miscov_{t-1}<-L_t+\alpha$.
Thus some $0\le r<t$ satisfies
\[
 S_{t-1}-S_r<-L_t+\alpha+1.
\]
For all sufficiently large $t$, $L_t\ge2(\alpha+1)$, so the right-hand side is at most $-L_t/2$.
The choice $r=t-1$ is then impossible.
For each fixed $r\le t-2$, the increment is the sum of $t-1-r$ independent, centered variables, each taking values in an interval of length $1$.
Then, Hoeffding's inequality \citep[Theorem~2, Eq.~(2.6), p.~16]{hoeffding1963probability}, applied to the negative increment, and a union bound give
\begin{align*}
 \PP\{t\text{ is forced from below}\}
 &\le\sum_{r=0}^{t-2}\exp\!\left(-\frac{L_t^2}{2(t-1-r)}\right)\\
 &\le t\exp\!\left(-\frac{c_-^2}{2}t^{1-2d}\right).
\end{align*}
Since $1-2d>0$, the last expression is summable in $t$.
Summing the trigger probabilities, including the finitely many initial rounds, proves $\EE N_\infty^-<\infty$.

For the other side, replace $e_t$ by $1-e_t$, $\alpha$ by $1-\alpha$, and exchange $(U_t,L_t)$.
The two triggers are mutually exclusive by \Cref{upper:stability}.
Induction in \eqref{eq:asym-truncation} therefore shows that the transformed output is $1-\hat e_t$ and its centered sum is $-\miscov_t$.
Its lower forced rounds are the original upper forced rounds.
The result already proved gives $\EE N_\infty^+<\infty$ when $b<1/2$.
\end{proof}

\subsection{Proof of Lemma~\ref{upper:mixed}}

\begin{proof}
Assume $d<1/2\le b$.
We bound the upper changes by the maximum of the free walk and the changes in the opposite direction.
If $M_T^+>0$, let $t\le T$ be its last upper change.
Then $M_T^+=M_{t-1}^++1$, while the upper trigger gives $\miscov_{t-1}>U_t-(1-\alpha)\ge-1$.
Using \eqref{upper:corrections}, we obtain
\[
 M_T^+=S_{t-1}+M_{t-1}^--\miscov_{t-1}+1
 <S_{t-1}+M_T^-+2.
\]
The case $M_T^+=0$ is immediate, so in all cases
\begin{equation}\label{upper:one-sided-count}
 M_T^+\le\max_{0\le r\le T}|S_r|+M_T^-+2.
\end{equation}
The free walk is a square-integrable martingale with $\EE S_T^2=\sigma^2T$.
Doob's $L^2$ maximal inequality \citep[Theorem~4.4.4, p.~236]{durrett2019probability} and Cauchy--Schwarz give
\[
 \EE\max_{0\le r\le T}|S_r|
 \le\left(\EE\max_{0\le r\le T}|S_r|^2\right)^{1/2}
 \le2\sigma\sqrt T.
\]
Also, $M_T^-\le N_\infty^-$ and \Cref{upper:loose} gives $\EE N_\infty^-<\infty$.
Taking expectations in \eqref{upper:one-sided-count} and using \eqref{upper:expected-changes},
\[
 \EE N_T^+
 \le\frac{2\sigma\sqrt T+\EE N_\infty^-+2}{\alpha}
 \le C\sqrt T\qquad(T\ge1).
\]
Reversing the upper and lower boundaries, as in the proof of \Cref{upper:loose}, gives the other case.
\end{proof}

\subsection{Proof of Lemma~\ref{upper:tight}}

\begin{proof}
If $\bunder=1$, then $N_T\le T$ proves the claim.
Assume $1/2\le\bunder<1$.
\paragraph{Step 1: bound the width and the midpoint's movement.}
Define
\[
 m_t=\frac{U_t-L_t}{2},\qquad
 w_t=\frac{U_t+L_t}{2},\qquad
 H_t=\miscov_t-m_t,\qquad
 \rho_t=m_t-m_{t-1}.
\]
Set $m_0=w_0=H_0=0$.
By feasibility, $|H_t|\le w_t$.
There are constants $0<c\le C<\infty$ such that
\begin{equation}\label{upper:potential-scales}
 c t^{1-\bunder}\le w_t\le C t^{1-\bunder},\qquad
 |\rho_t|\le C t^{-\bunder},\qquad
 w_{t-1}|\rho_t|\le C
 \qquad(t\ge1).
\end{equation}
Indeed, one of $b,d$ equals $\bunder$, so the first bound holds with lower constant $\min\{c_-,c_+\}/2$ and upper constant $(c_-+c_+)/2$.
For $t\ge2$ and $\gamma\in[0,1]$, the mean value theorem gives 
\[0\le t^\gamma-(t-1)^\gamma\le2t^{\gamma-1}.\]
Applying this to $1-b$ and $1-d$ proves the second bound, provided the constant is increased to cover $t=1$.
The first two bounds give $w_{t-1}|\rho_t|\le C t^{1-2\bunder}\le C$.
Also, $w_t$ is nondecreasing and tends to infinity, while $\rho_t$ tends to zero.
Choose $t_\star$ so that
\begin{equation}\label{upper:large-width}
 w_t\ge2,\qquad |\rho_t|\le\nu/2
 \qquad(t\ge t_\star).
\end{equation}

\paragraph{Step 2: bound the increase on unforced rounds.}
Let $I_t^0$ indicate that neither trigger holds.
On this event,
\[
 H_t=H_{t-1}+e_t-\alpha-\rho_t.
\]
Both $I_t^0$ and $H_{t-1}$ are determined by the previous input errors, whereas $e_t-\alpha$ is independent of them, centered, and has variance $\sigma^2$.
Thus
\begin{align}
 \EE\bigl[(H_t^2-H_{t-1}^2)I_t^0\bigr]
 &=\EE\bigl[(-2H_{t-1}\rho_t+\sigma^2+\rho_t^2)I_t^0\bigr]
 \le C_1,\label{upper:unforced-drift}
\end{align}
where the last bound follows from $|H_{t-1}|\le w_{t-1}$ and \eqref{upper:potential-scales}.

\paragraph{Step 3: bound the decrease on forced rounds.}
Fix $t\ge t_\star$ and put $h=\miscov_{t-1}-m_t$, so $H_{t-1}=h+\rho_t$.
If the upper trigger holds, then $h>w_t-(1-\alpha)$ and $H_t=h-\alpha$.
Consequently,
\[
 H_t^2-H_{t-1}^2
 =-(\alpha+\rho_t)(2h-\alpha+\rho_t)
 \le-\frac{\alpha}{2}w_t.
\]
Here \eqref{upper:large-width} gives $\alpha+\rho_t\ge\alpha/2$ and $2h-\alpha+\rho_t\ge2w_t-2\ge w_t$.
If the lower trigger holds, then $h<-w_t+\alpha$ and $H_t=h+1-\alpha$.
In this case,
\[
 H_t^2-H_{t-1}^2
 =(1-\alpha-\rho_t)(2h+1-\alpha+\rho_t)
 \le-\frac{1-\alpha}{2}w_t,
\]
because the first factor is at least $(1-\alpha)/2$ and the second is at most $-2w_t+2\le-w_t$.
Writing $I_t=1-I_t^0$ for the indicator of a forced round, both cases give
\begin{equation}\label{upper:forced-drift}
 (H_t^2-H_{t-1}^2)I_t\le-\frac{\nu}{2}w_t I_t
 \qquad(t\ge t_\star).
\end{equation}

\paragraph{Step 4: count the forced rounds in a block.}
For an integer $n\ge t_\star$, sum \eqref{upper:unforced-drift} and the expectation of \eqref{upper:forced-drift} over $n<t\le2n$.
Since $w_t\ge w_n$ on the block,
\[
 \EE H_{2n}^2-\EE H_n^2
 \le C_1n-\frac{\nu}{2}w_n\EE(N_{2n}-N_n).
\]
The left side is at least $-w_n^2$.
Rearranging and using \eqref{upper:potential-scales},
\[
 \EE(N_{2n}-N_n)
 \le\frac{2}{\nu}\left(\frac{C_1n}{w_n}+w_n\right)
 \le C\bigl(n^{\bunder}+n^{1-\bunder}\bigr)
 \le C n^{\bunder}.
\]
The last inequality uses $\bunder\ge1/2$.

\paragraph{Step 5: sum the blocks.}
Choose $j_0\ge0$ with $2^{j_0}\ge t_\star$.
For $i>j_0$, Step~4 and $N_{2^{j_0}}\le2^{j_0}$ give
\[
 \EE N_{2^i}
 \le2^{j_0}+C\sum_{j=j_0}^{i-1}2^{j\bunder}
 \le C'2^{i\bunder}.
\]
After increasing $C'$, this holds also for $0\le i\le j_0$.
For any $T\ge1$, choose $i\ge0$ with $2^{i-1}<T\le2^i$.
Monotonicity of $N_T$ then gives $\EE N_T\le\EE N_{2^i}\le C'(2T)^{\bunder}$, which proves the claim.
\end{proof}

\section{Unconstrained minimax rate}
\label{learning:section}

We use the notation from Appendix~\ref{app:asymmetric-setting}. We first consider the following lemmas.
For a fixed distribution, the first lemma turns disagreement with the closest exact oracle into a sum of instantaneous training-conditional errors.

\begin{lemma}
\label{learning:repair}
Fix $\PP\in\calQ$ and a deterministic measurable online algorithm $\calC$.
There is a deterministic measurable oracle $\calC^{\mathrm{rep}}_{\PP}\in\calO_\alpha(\PP)$ such that, for every $T\ge1$,
\[
 r_T(\calC;\PP)
 =\EE_{\PP}\ell_T(\calC,\calC^{\mathrm{rep}}_{\PP})
 =\sum_{t=1}^T\EE_{\PP}|p_t^{\calC}-\alpha|.
\]
The same oracle works for every horizon.
\end{lemma}

Let $B_t\sim\operatorname{Binomial}(t,\alpha)$.
For $t$ i.i.d.\ observations, $B_t$ has the distribution of the number of miscoverages by the population score oracle.
If the realized count is $b$, declaring miscoverage on the $k$ largest scores disagrees with the population oracle on exactly $|b-k|$ observations.
The following lemma shows why the same count yields a lower bound even for algorithms that need not threshold scores.

\begin{lemma}[Comparison with a binomial count]\label{learning:binomial-lower}
For every $\calC\in\calA$ and every $T\ge1$,
\[
 R_T(\calC)\ge\sum_{t=1}^T\frac1t
       \min_{k\in\{0,\ldots,t\}}\EE|B_t-k|.
\]
\end{lemma}

The last lemma supplies the upper bound using the unsmoothed conformal set in the main text.

\begin{lemma}[Conformal disagreement with the population oracle]
\label{learning:conformal-upper}
For every $\PP\in\calQ$ and every $t\ge1$,
\[
 \EE_{\PP}\left|
 \PP\{Y_t\in\calC_t^{\cp}(X_t)\mid\calG_{t-1}\}-(1-\alpha)
 \right|
 \le \frac{1}{2\sqrt{t+1}}+\frac1t.
\]
The population score oracle attains the smallest expected disagreement with $\calC^{\cp}$ among all exact oracles.
Consequently, for every $T\ge1$,
\[
 R_T(\calC^{\cp})\le3\sqrt T.
\]
\end{lemma}

\subsection{Proof of Lemma~\ref{lem:learning-full}}
\label{learning:main-proof}

\begin{proof}
We first match the binomial lower bound by a prediction set with an integer cutoff at each round.
The rate follows from the size of the binomial fluctuations.

\paragraph{Step 1. Match the lower bound at every horizon.}
Fix deterministic integers $0\le k_t\le t$ and form the prediction set
\[
 \calC^{(k)}_t(x)=
 \left\{y:1+\sum_{r<t}\ind\{s_r\ge s(x,y)\}>k_t\right\}.
\]
Under any $\PP\in\calQ$, the $t$ scores are almost surely distinct.
Hold out each observation in turn and use the other $t-1$ scores in this formula.
The algorithm declares miscoverage on the $k_t$ largest scores.
The population oracle declares miscoverage on the scores above $q_{\PP,s}$, whose number has distribution $B_t$.
The two subsets are nested.
By exchangeability, the disagreement probability at round $t$ equals the expected fraction of these $t$ decisions that disagree.
Thus, for every $\PP\in\calQ$,
\[
 \EE_\PP\ell_T(\calC^{(k)},\calC^*_{\PP,s})
 =\sum_{t=1}^T\frac{\EE|B_t-k_t|}{t}.
\]

Choose an integer median of $B_t$ for $k_t$.
It minimizes the summand because $\EE|B_t-(k+1)|-\EE|B_t-k|=2\Pr(B_t\le k)-1$.
Comparing with the population oracle gives an upper bound uniform over $\calQ$.
Together with \Cref{learning:binomial-lower}, this gives the identity
\begin{equation}\label{learning:binomial-minimax}
 R_T^0=\sum_{t=1}^T\frac1t\min_{k\in\{0,\ldots,t\}}\EE|B_t-k| .
\end{equation}

\paragraph{Step 2. Bound the binomial fluctuations.}
Put $X=B_t-t\alpha$ and let $\xi\sim\operatorname{Bernoulli}(\alpha)$.
Independence and centering give
\[
 \EE X^2=t\sigma^2,\qquad
 \EE X^4=t\EE(\xi-\alpha)^4+3t(t-1)\sigma^4\le 2t^2\sigma^2,
\]
using $|\xi-\alpha|\le1$ and $\sigma^2\le1/4$.
For any $k$, $|t\alpha-k|\le\EE|B_t-k|$ by Jensen's inequality.
The triangle inequality and H\"older's inequality therefore give
\[
 \EE|B_t-k|\ge\tfrac12\EE|X|
 \ge\frac{(\EE X^2)^{3/2}}{2(\EE X^4)^{1/2}}
 \ge\frac{\sigma^2}{2\sqrt2}\sqrt t .
\]
Since $\sigma^2\ge\nu^2$ and $\sum_{t=1}^Tt^{-1/2}\ge\sqrt T$, \eqref{learning:binomial-minimax} implies $R_T^0\ge(\nu^2/32)\sqrt T$.

The usual split conformal set in \eqref{eq:conformal-baseline-set} corresponds to $k_t=\lfloor\alpha t\rfloor$.
For this choice, $\EE|B_t-k_t|\le\sigma\sqrt t+1$.
Since $\sigma\le1/2$ and $\sum_{t=1}^T t^{-1/2}\le2\sqrt T$, Step~1 gives $R_T^0\le R_T(\calC^{\cp})\le3\sqrt T$.
\end{proof}

\subsection{Proof of Lemma~\ref{learning:repair}}
\label{learning:repair-proof}

\begin{proof}
We present the lemma in the following steps.

\paragraph{Step 1. Bound the disagreement of any exact oracle.}
For the error $e_t^*$ of any exact oracle, conditional Jensen's inequality gives
\[
 \EE_\PP[|\err_t^{\calC}-e_t^*|\mid\calG_{t-1}]
 \ge\bigl|\EE_\PP[\err_t^{\calC}-e_t^*\mid\calG_{t-1}]\bigr|
 =|p_t^{\calC}-\alpha|.
\]
It suffices to attain this bound at each round with the same oracle.

\paragraph{Step 2. Construct the oracle and compute its loss.}
Fix a round $t$ and history $h$.
For a fresh $Z\sim\mu_\PP$, put
\[
 A_h(z)=e_t(h,z)+F(s(z)),\qquad
 G_h(a)=\int\ind\{A_h(z')\le a\}\,d\mu_\PP(z').
\]
Since $F(s(Z))$ is uniform on $(0,1)$, the algorithm covers exactly when $A_h(Z)\le1$, almost surely.
Also $\{A_h(Z)=a\}\subseteq\{F(s(Z))\in\{a,a-1\}\}$, so $A_h(Z)$ has no atoms and $G_h$ is continuous.

Define
\[
 \calC^{\mathrm{rep}}_{\PP,t}(x;h)
 =\{y:G_h(A_h(x,y))\le1-\alpha\}.
\]
Given $h$, the probability integral transform makes $G_h(A_h(Z))$ uniformly distributed, so the oracle has instantaneous training-conditional coverage of $1-\alpha$.

Both coverage events are lower intervals in $A_h(Z)$ and hence are nested.
If $p(h)=\int e_t(h,z)\,d\mu_\PP(z)$, their probabilities are $G_h(1)=1-p(h)$ and $1-\alpha$.
Their conditional disagreement is therefore
\[
 \bigl|G_h(1)-(1-\alpha)\bigr|=|p(h)-\alpha|.
\]
At the observed history $p(h)=p_t^{\calC}$ almost surely.
Taking expectations and summing attains the lower bound from Step~1 for every $T$, since the construction is independent of the horizon.
\end{proof}

\subsection{Proof of Lemma~\ref{learning:binomial-lower}}
\label{learning:binomial-proof}

\begin{proof}
Fix $\calC$ and $T$. We first compute the calibration error of a decision that depends only on the relative order of the scores. We then choose one distribution under which such decisions approximate $\calC$ at every $t\le T$.

\paragraph{Step 1. Count disagreements for a rule based on ranks.}
Let $V_1,\ldots,V_t$ be i.i.d.\ uniform on $(0,1)$. Write $R_i$ for the rank of $V_i$ among all $t$ values, and let $B_t=\sum_{i=1}^t\ind\{V_i\le\alpha\}$.
For $0\le k\le t$, the two indicators below select the first $k$ and the first $B_t$ ranks. Hence
\[
 \sum_{i=1}^t\bigl|\ind\{R_i\le k\}-\ind\{V_i\le\alpha\}\bigr|=|B_t-k|.
\]
Let $V_{(1)},\ldots,V_{(t-1)}$ be the sorted past values, with $V_{(0)}=0$ and $V_{(t)}=1$.
Since $R_t\le k$ exactly when $V_t\le V_{(k)}$ almost surely, conditioning on the past and then using exchangeability gives
\[
 \EE|V_{(k)}-\alpha|
 =\EE\bigl|\ind\{R_t\le k\}-\ind\{V_t\le\alpha\}\bigr|
 =\frac1t\EE|B_t-k|.
\]

Fix the past order $\rho$. A rank-based rule declares miscoverage on some set $J(\rho)$. Given the past values, its miscoverage probability is
\[
 q_t=\sum_{j\in J(\rho)}G_j,\qquad G_j=V_{(j)}-V_{(j-1)} .
\]
Given $\rho$, the sorted values have constant density $(t-1)!$, so the gaps are uniform on the simplex $\{g_j\ge0,\ \sum_jg_j=1\}$ and are exchangeable. For $t=1$ the sole gap is $G_1=1$.
For $k=|J(\rho)|$, the sum has the same distribution as $\sum_{j=1}^kG_j=V_{(k)}$. Thus for each $\rho$,
\[
 \EE[|q_t-\alpha|\mid\rho]
 =\frac1t\EE|B_t-k|
 \ge \frac1t\min_{0\le j\le t}\EE|B_t-j|.
\]

\paragraph{Step 2. Approximate $\calC$ using only score ranks.}
Fix an integer $m\ge T$. We will give equal mass to $m$ cells ordered by score. On each tuple of distinct cells, a function of the score ranks will approximate the miscoverage indicator $e_t=\ind\{Y_t\notin\calC_t(X_t)\}$ with error probability below $1/m$.

\emph{Choose the cells.}
Let $N$ be large enough. 
Take $\mu_0$ with $\mu_0^{\otimes\mathbb N}\in\calQ$. Since the score is nonatomic, we can split $\calZ$ into $N$ consecutive score bins $I_1,\ldots,I_N$, each of mass $1/N$.
For each $t\le T$, approximate $\{e_t=1\}$ in product measure by a finite union of rectangles with indicator $f_t$. Refine the bins by the sides of all these rectangles. Then $f_t$ is constant on each product cell, and we can make the errors $E_t=\{e_t\ne f_t\}$ small enough that
\[
 \sum_{t=1}^T N^t\mu_0^{\otimes t}(E_t)<1/m .
\]

Independently choose one partition cell $D_i$ in each bin $I_i$, with $\Pr(D_i=D)=N\mu_0(D)$ for each cell $D\subseteq I_i$. Set $\mu_i=\mu_0(\cdot\mid D_i)$. Averaging $\mu_i$ over this choice gives $\mu_0(\cdot\mid I_i)$.
For each ordered tuple $\mathbf i=(i_1,\ldots,i_t)$ of distinct labels, write $\mu_{\mathbf i}=\bigotimes_{r=1}^t\mu_{i_r}$ and $I_{\mathbf i}=\prod_{r=1}^t I_{i_r}$. Independence and disjointness of the bin products at each $t$ give
\[
 \EE\sum_{t=1}^T\sum_{\mathbf i}\mu_{\mathbf i}(E_t)
 =\sum_{t=1}^T N^t\sum_{\mathbf i}\mu_0^{\otimes t}(E_t\cap I_{\mathbf i})
 \le\sum_{t=1}^T N^t\mu_0^{\otimes t}(E_t)<1/m .
\]
Fix a choice of cells for which the sum is below $1/m$. Then $\mu_{\mathbf i}(E_t)<1/m$ for every tuple in the sum.

\emph{Select cells on which ranks suffice.}
The color of a set of $t$ labels records which orders give $f_t=1$. For $i_1<\cdots<i_t$, define
\[
 c_t(i_1,\ldots,i_t)
 =\bigl\{\pi\in S_t:\ f_t=1\text{ on }
       D_{i_{\pi(1)}}\times\cdots\times D_{i_{\pi(t)}}\bigr\},
\]
where $S_t$ is the set of permutations of $\{1,\ldots,t\}$.
There are at most $2^{t!}$ colors. Applying the finite Ramsey theorem \citep[Theorem~B, p.~267]{ramsey1930problem} successively for $t=1,\ldots,T$ gives a set $H$ of $m$ labels such that, for each $t$, every $t$-subset of $H$ has the same color $\Pi_t$.

For distinct scores, let $\pi(r)$ be the rank of $s_r$ among $s_1,\ldots,s_t$ and set $\tilde e_t=\ind\{\pi\in\Pi_t\}$. On distinct cells from $H$, the scores have the same relative order as their cell labels. Thus $f_t=\tilde e_t$ on each such product cell, and
\[
 \mu_{\mathbf i}\{e_t\ne\tilde e_t\}
 =\mu_{\mathbf i}(E_t)<1/m
 \qquad (i_1,\ldots,i_t\text{ distinct in }H).
\]

\paragraph{Step 3. Sum the bound under one distribution.}
Relabel the retained cells $1,\ldots,m$ and set $\mu=m^{-1}\sum_{i=1}^m\mu_i$. Since $\mu\ll\mu_0$, the score remains nonatomic and $\PP=\mu^{\otimes\mathbb N}\in\calQ$.
Set $\tilde p_t=\EE_\PP[\tilde e_t\mid\calG_{t-1}]$. The cell labels are i.i.d.\ uniform on $\{1,\ldots,m\}$, so Step~2 and a union bound give
\[
 \EE_\PP|p_t^{\calC}-\tilde p_t|
 \le \PP(\err_t^{\calC}\ne\tilde e_t)
 \le \frac1m+\PP(\text{two of the first $t$ labels coincide})
 \le \frac{1+\binom t2}{m}.
\]
The uniforms $V_r=1-F(s_r)$ rank the scores in decreasing order, and $V_t$ is independent of $\calG_{t-1}$. Thus $\tilde p_t$ is the gap sum in Step~1. By \Cref{learning:repair} and the triangle inequality,
\[
 R_T(\calC)\ge r_T(\calC;\PP)
 =\sum_{t=1}^T\EE_\PP|p_t^{\calC}-\alpha|
 \ge \sum_{t=1}^T\frac1t\min_{0\le k\le t}\EE|B_t-k|
     -\frac1m\sum_{t=1}^T\left(1+\binom t2\right).
\]
Let $m\to\infty$ with $T$ fixed.
\end{proof}

\subsection{Proof of Lemma~\ref{learning:conformal-upper}}
\label{learning:conformal-proof}

\begin{proof}
Fix $\PP\in\calQ$.
Let $k_t=\lceil(1-\alpha)t\rceil$.
If $k_t\le t-1$, the unsmoothed conformal set in \eqref{eq:conformal-baseline-set} is, up to a fresh-score null event, the set whose cutoff is the $k_t$-th smallest of the $t-1$ past scores.

The transformed scores $F(s(Z_r))$ are independent uniform random variables.
Consequently, the instantaneous training-conditional coverage probability $Q_t=1-p_t^{\calC^{\cp}}$ of the conformal set has distribution $\operatorname{Beta}(k_t,t-k_t)$ when $k_t<t$, and equals one when $k_t=t$.
The beta mean and variance, together with the case $k_t=t$, give
\[
 \EE Q_t=\frac{k_t}{t},
 \qquad
 \operatorname{Var}(Q_t)\le\frac1{4(t+1)},
 \qquad
 \left|\frac{k_t}{t}-(1-\alpha)\right|\le\frac1t.
\]
The triangle inequality and Cauchy--Schwarz prove
\[
 \EE|Q_t-(1-\alpha)|
 \le\frac1{2\sqrt{t+1}}+\frac1t.
\]
The conformal and population prediction sets are thresholds of the same score, with both cutoffs determined before the fresh outcome.
Their conditional disagreement probability is therefore exactly $|Q_t-(1-\alpha)|$.
Conditional Jensen's inequality gives this same quantity as a lower bound for every exact oracle, so the population oracle attains the infimum.
Finally,
\[
 \sum_{t=1}^T\left(\frac1{2\sqrt{t+1}}+\frac1t\right)
 \le\frac32\sum_{t=1}^Tt^{-1/2}\le3\sqrt T.
\]
The bound is uniform over $\PP\in\calQ$, which proves the lemma.
\end{proof}

\section{Results with a relaxed oracle condition}
\label{app:relaxed}

An exact oracle condition may be impossible when the joint distribution has atoms.
We therefore also consider relaxed oracles whose instantaneous training-conditional coverage is at least $1-\alpha$:
\begin{equation}\label{eq:relaxed-oracle}
 \PP\{Y_t\in\calC^*_{\PP,t}(X_t)\mid\calG_{t-1}\}
 \ge1-\alpha
 \quad\text{almost surely for every }t\ge1.
 \tag{oracle$'$}
\end{equation}
Let $\calO_\alpha^+(\PP)$ be a nonempty family of such oracles, fixed independently of the algorithm and horizon.
The population quantile set \eqref{eq:population-oracle} satisfies this condition even when the score distribution has atoms.
For a specified nonempty family $\calO_\alpha^+(\PP)$ satisfying \eqref{eq:relaxed-oracle}, write
\[
 R_T^+(\calC)=\sup_{\PP\in\tilde{\calQ}}\inf_{\calC^*_{\PP}\in\calO_\alpha^+(\PP)}
                 \EE_\PP\ell_T(\calC,\calC^*_{\PP}),\qquad
 R_T^{0,+}=\inf_{\calC\in\calA}R_T^+(\calC).
\]

For this relaxed formulation, every $\calC \in \Aasym$ satisfies the lower-bound decomposition over the class $\calQ$ of all i.i.d. laws:
\begin{equation}\label{eq:relaxed-lower-decomposition}
\begin{split}
 &\sup_{\PP\in\tilde\calQ}
 \inf_{\calC^*_{\PP}\in\calO_\alpha^+(\PP)}
 \EE_{\PP}\ell_T(\calC,\calC^*_{\PP})\\
 &\qquad\gtrsim
 \inf_{\calC^0\in\calA}\sup_{\PP\in\tilde\calQ}
 \inf_{\calC^*_{\PP}\in\calO_\alpha^+(\PP)}
 \EE_{\PP}\ell_T(\calC^0,\calC^*_{\PP})
 +
 \begin{cases}
  0,&0\le\beta<1/2,\\
  T^\beta,&1/2\le\beta\le1.
 \end{cases}
\end{split}
\end{equation}
The unconstrained term need not have order $\sqrt T$.
For example, if $\calO_\alpha^+(\PP)$ contains all relaxed oracles, this term is zero: the algorithm that always predicts $\calY$ agrees perfectly with the oracle that does the same.
A relaxed oracle may therefore be a poor target for statistical accuracy, since matching it need not yield informative predictions or coverage close to $1-\alpha$.

Formally, we have the following result.
\begin{proposition}[Lower bounds with relaxed oracles]\label{prop:asym-relaxed}
Every $\calC\in\Aasym$ satisfies, for every $T\ge1$,
\begin{equation}\label{eq:asym-relaxed-lower}
 R_T^+(\calC)\ge\tfrac12\left(R_T^{0,+}+(\alpha T-L_T)_+\right).
\end{equation}
If $d>0$, or if $d=0$ and $c_-<\alpha$, the constrained minimax risk has order $T$.
If the oracle family contains all relaxed oracles and $d=0$, $c_-\ge\alpha$, the constrained minimax risk is zero.
\end{proposition}
Only the lower boundary forces a prediction set to miss a response that an always-covering oracle covers.
This explains why $g_{b,d}$ is not the appropriate universal lower bound for relaxed oracles.

To prove the lower bound above, we first consider the following lemma.
\begin{lemma}[Disagreement under a point mass]\label{relaxed:point-mass}
Fix any $z_0=(x_0,y_0)\in\calX\times\calY$ and let $\PP_{\rm pt}$ be the i.i.d.\ law concentrated on $z_0$.
For every $\calC\in\Aasym$, every deterministic relaxed oracle $\calC^*_{\PP}$ under $\PP_{\rm pt}$, and every $T\ge1$,
\[
 \EE_{\PP_{\rm pt}}\ell_T(\calC,\calC^*_{\PP})
 =\sum_{t=1}^T\err_t^{\calC}
 \ge(\alpha T-L_T)_+.
\]
\end{lemma}

\subsection{Proof of Proposition~\ref{prop:asym-relaxed}}
\label{relaxed:proof}

\begin{proof}
We prove the proposition in the following steps.

\paragraph{Step 1. Combine the two lower bounds.}
Fix $\calC\in\Aasym$ and $T\ge1$.
Since $\calC\in\calA$, we have $R_T^+(\calC)\ge R_T^{0,+}$.
Since $\PP_{\rm pt}\in\tilde{\calQ}$, \Cref{relaxed:point-mass} gives
\[
 R_T^+(\calC)\ge
 \inf_{\calC^*_{\PP}\in\calO_\alpha^+(\PP_{\rm pt})}
     \EE_{\PP_{\rm pt}}\ell_T(\calC,\calC^*_{\PP})
 \ge(\alpha T-L_T)_+ .
\]
Taking the larger of the bounds yields
\[
 R_T^+(\calC)\ge\max\{R_T^{0,+},(\alpha T-L_T)_+\}
 \ge\tfrac12\bigl(R_T^{0,+}+(\alpha T-L_T)_+\bigr),
\]
which proves \eqref{eq:asym-relaxed-lower} for every horizon.

\paragraph{Step 2. Determine when the forced cost is linear.}
Since $L_T=c_-T^{1-d}$,
\[
 \frac{(\alpha T-L_T)_+}{T}=(\alpha-c_-T^{-d})_+
 \longrightarrow
 \begin{cases}
  \alpha,&d>0,\\
  (\alpha-c_-)_+,&d=0 .
 \end{cases}
\]
The limit is positive if $d>0$ or if $d=0$ and $c_-<\alpha$.
Truncation supplies a feasible algorithm by \Cref{upper:stability}, and $\ell_T\le T$ on every path.
Hence in either case, for all sufficiently large $T$,
\[
 cT\le\inf_{\calC\in\Aasym}R_T^+(\calC)\le T
\]
for some $c>0$.

\paragraph{Step 3. Attain zero risk when the lower boundary permits full coverage.}
Suppose $d=0$, $c_-\ge\alpha$, and the comparison family contains all relaxed oracles.
The algorithm $\calC^{\rm all}_t(x)=\calY$ satisfies
\[
 -L_t=-c_-t\le\miscov_t^{\calC^{\rm all}}=-\alpha t\le U_t
 \qquad(t\ge1),
\]
so it belongs to $\Aasym$.
Its instantaneous training-conditional coverage is one under every distribution, so it is also an admissible comparison oracle for each $\PP\in\tilde{\calQ}$.
Comparing the algorithm with itself gives
\[
 0\le\inf_{\calC\in\Aasym}R_T^+(\calC)\le R_T^+(\calC^{\rm all})
 \le\sup_{\PP\in\tilde{\calQ}}\EE_\PP\ell_T(\calC^{\rm all},\calC^{\rm all})=0 .
\]
\end{proof}

\subsection{Proof of Lemma~\ref{relaxed:point-mass}}
\label{relaxed:point-mass-proof}

\begin{proof}
On the constant sequence $Z_t=z_0$, both the history $h=(z_0,\ldots,z_0)$ through round $t-1$ and the fresh response $y_0$ are fixed.
For a deterministic relaxed oracle, \eqref{eq:relaxed-oracle} becomes
\[
 1-\alpha
 \le \PP_{\rm pt}\{Y_t\in\calC_t^*(X_t)\mid\calG_{t-1}\}
 =\ind\{y_0\in\calC_t^*(x_0;h)\}\in\{0,1\}.
\]
Since $1-\alpha>0$, the indicator is one at every round.
The oracle's error sequence is therefore zero, and on the unique data path
\[
 \EE_{\PP_{\rm pt}}\ell_T(\calC,\calC^*_{\PP})
 =\sum_{t=1}^T|\err_t^{\calC}-0|
 =\alpha T+\miscov_T^{\calC}
 \ge\alpha T-L_T .
\]
The error count is also nonnegative, giving the positive part in the claim.
\end{proof}

Take $b=d=\beta$ and $c_-=c_+=\kappa$.
If $\beta\ge1/2$, then $\beta>0$, so Step~2 gives $R_T^+(\calC)\gtrsim T\ge T^\beta$ for every constrained algorithm.
Together with $R_T^+(\calC)\ge R_T^{0,+}$, this proves the second case of \eqref{eq:relaxed-lower-decomposition}.
If $\beta<1/2$, that display adds zero to the unconstrained risk, so class inclusion alone proves its first case.
In fact, the point-mass argument gives the stronger linear lower bound whenever $\beta>0$.
For asymmetric constraints, a tight upper bound alone does not yield this lower bound: Step~3 shows that zero risk is possible if the lower bound permits the always-covering algorithm.

\section{Proofs for ACI}
\label{app:aci}\label{sec:proof-aci}

This section proves \Cref{prop:aci-subopt,cor:aci-decaying-risk}.
Here $s_t\in[0,B]$ and $\err_t=\ind\{s_t>q_t\}$ is the ACI error.
The update \eqref{eq:aci1} is $q_{t+1}=q_t+\eta_t(\err_t-\alpha)$, where $q_1\in[0,B]$ and the positive step sizes are deterministic.
Write $\bar\eta_t=\max_{1\le r\le t}\eta_r$, with $\bar\eta_0=0$.
The following two bounds hold for every score sequence, without any distributional assumption.
We cite their proofs from \citet{angelopoulos2024online}.

\begin{lemma}[{\citealp[Lemma~1]{angelopoulos2024online}}]\label{lem:aci_range}
For every sequence of scores in $[0,B]$ and every $t\ge1$,
\[
 -\alpha\bar\eta_{t-1}\le q_t\le B+(1-\alpha)\bar\eta_{t-1}.
\]
\end{lemma}

\begin{lemma}[{\citealp[Theorems~1 and~2]{angelopoulos2024online}}]\label{lem:aci_corridor}
For every sequence of scores in $[0,B]$ and every $T\ge1$,
\[
 |\miscov_T^{\calC^{\aci}}|
 \le (B+\bar\eta_T)\left(\frac1{\eta_1}
      +\sum_{t=2}^T\left|\frac1{\eta_t}-\frac1{\eta_{t-1}}\right|\right).
\]
For nonincreasing step sizes, this becomes $|\miscov_T^{\calC^{\aci}}|\le(B+\eta_1)/\eta_T$.
In particular, $\eta_t=ct^{-(1-\beta)}$ with $c>0$ and $\beta\in(0,1]$ gives the $(\beta,\kappa)$ adversarial guarantee with $\kappa=(B+c)/c$.
\end{lemma}

The cited theorems bound the long-run empirical coverage rate, with the difference from $1-\alpha$ given by $-\miscov_T^{\calC^{\aci}}/T$.
For nonincreasing steps, $1/\eta_1+\sum_{t=2}^T|1/\eta_t-1/\eta_{t-1}|=1/\eta_T$ and $\bar\eta_T=\eta_1$.
This gives \eqref{eq:aci-worst-case} for the stated power steps.

\paragraph{The score distribution and step sizes.}
For the rest of this section, fix an i.i.d.\ law $\PP$ under which the score has a density satisfying \eqref{eq:aci-density}, and suppose
\begin{equation}\label{eq:aci-step-assumptions}
 \eta_t\longrightarrow0,
 \qquad
 \left|\eta_{t+1}^{-1}-\eta_t^{-1}\right|\longrightarrow0.
\end{equation}
In particular, $\bar\eta:=\sup_{t\ge1}\eta_t<\infty$.
The score CDF $F$ equals $0$ below $0$ and $1$ above $B$.
The density bounds give
\begin{equation}\label{eq:aci_lip}
 \lambda_-(u-v)\le F(u)-F(v)\le\lambda_+(u-v)
 \qquad(0\le v\le u\le B).
\end{equation}
Thus $F$ is continuous and strictly increasing on $[0,B]$, so $q_{\PP,s}\in(0,B)$ and $F(q_{\PP,s})=1-\alpha$.

\paragraph{Cumulative quantile deviation and instantaneous training-conditional miscoverage.}
The threshold $q_t$ is determined by the past observations, and the fresh score $s_t$ is independent of that past.
Hence, almost surely,
\begin{equation}\label{eq:aci_pt}
 p_t=p_t^{\calC^{\aci}}=1-F(q_t).
\end{equation}
Define
\begin{equation}\label{eq:aci_notation}
 \Delta_t=q_t-q_{\PP,s},
 \qquad
 \chi(u)=F(q_{\PP,s}+u)-F(q_{\PP,s}).
\end{equation}
Then $p_t-\alpha=-\chi(\Delta_t)$, and the recursion becomes
\begin{equation}\label{eq:aci_delta}
 \Delta_{t+1}=\Delta_t-\eta_t\chi(\Delta_t)+\eta_t\xi_t,
 \qquad \xi_t=\err_t-p_t.
\end{equation}
Here $-\eta_t\chi(\Delta_t)$ pulls the threshold toward the population quantile, while $\eta_t\xi_t$ is noise satisfying
\begin{equation}\label{eq:aci_noise}
 \EE_{\PP}[\xi_t\mid\calG_{t-1}]=0,
 \qquad
 \EE_{\PP}[\xi_t^2\mid\calG_{t-1}]=p_t(1-p_t),
 \qquad |\xi_t|\le1.
\end{equation}
We suppress the subscript $\PP$ on expectations below.

\paragraph{From disagreement to threshold deviations.}
The next lemma identifies the fixed-distribution risk in \eqref{eq:asym-risk-notation} and shows that the population prediction set attains it.

\begin{lemma}[Disagreement with the closest exact oracle]\label{lem:aci_risk}
For every $T\ge1$,
\begin{equation}\label{eq:aci_risk}
 r_T(\calC^{\aci};\PP)
 =\EE\sum_{t=1}^T|p_t-\alpha|
 =\EE\sum_{t=1}^T|\chi(\Delta_t)|.
\end{equation}
The population prediction set in \eqref{eq:population-oracle} attains the infimum.
\end{lemma}

The density bounds make the change in instantaneous training-conditional miscoverage comparable to the threshold deviation, including when the threshold leaves $[0,B]$.

\begin{lemma}[Comparison with the threshold deviation]\label{lem:aci_lip}
Let $\underline\lambda=\min\{\lambda_-,\nu/(B+\bar\eta)\}>0$.
Then $|\Delta_t|\le B+\bar\eta$ for every $t\ge1$.
For every $u$ with $|u|\le B+\bar\eta$, the value $\chi(u)$ has the same sign as $u$ and
\[
 \underline\lambda|u|\le|\chi(u)|\le\lambda_+|u|.
\]
\end{lemma}

The last lemma gives the expected size of the deviation.
Its proof balances the contraction of the second moment, of order $\eta_t\EE\Delta_t^2$, against noise of order $\eta_t^2$.
A fourth-moment bound ensures that the second moment does not come only from rare large deviations.

\begin{lemma}[Expected absolute deviation]\label{lem:aci_fluctuation}
Under \eqref{eq:aci-step-assumptions}, there are constants $0<c_0\le C_0<\infty$ and a time $t_\star\ge1$ such that
\[
 c_0\sqrt{\eta_t}\le\EE|\Delta_t|\le C_0\sqrt{\eta_t}
 \qquad(t\ge t_\star).
\]
\end{lemma}

We first combine these lemmas to prove the proposition and derive the corollary.
The supporting lemma proofs follow.

\subsection{Proof of Proposition~\ref{prop:aci-subopt}}
\label{app:aci-subopt-proof}

\begin{proof}
The oracle identity and the density bounds reduce expected disagreement to the sum of expected absolute threshold deviations.
We apply the deviation estimate and show that finitely many initial rounds are negligible.
Write $A^{\eta}_T=\sum_{t=1}^T\sqrt{\eta_t}$, with $A^{\eta}_0=0$.

\paragraph{Step 1: from the risk to the deviations.}
By \Cref{lem:aci_risk} and \Cref{lem:aci_lip},
\begin{equation}\label{eq:aci_sandwich}
  \underline{\lambda}\sum_{t=1}^T\EE|\Delta_t|
  \;\le\; r_T(\calC^{\aci};\PP)\;\le\;
  \lambda_+\sum_{t=1}^T\EE|\Delta_t| .
\end{equation}

\paragraph{Step 2: the size of the deviations.}
By \Cref{lem:aci_fluctuation}, $c_0\sqrt{\eta_t}\le\EE|\Delta_t|\le C_0\sqrt{\eta_t}$ for $t\ge t_\star$.
For $t<t_\star$ we use $0\le\EE|\Delta_t|\le B+\bar\eta$ from \Cref{lem:aci_lip}.
Inserting these bounds into \eqref{eq:aci_sandwich} gives, for $T\ge t_\star$,
\[
  \underline{\lambda}\,c_0\bigl(A^{\eta}_T-A^{\eta}_{t_\star-1}\bigr)
  \;\le\; r_T(\calC^{\aci};\PP)\;\le\;
  \lambda_+C_0\,A^{\eta}_T+\lambda_+(B+\bar\eta)\,t_\star .
\]

\paragraph{Step 3: the first rounds are negligible.}
It remains to show that $A^{\eta}_T\to\infty$, so that the constants $A^{\eta}_{t_\star-1}$ and $\lambda_+(B+\bar\eta)t_\star$ do not matter.
Condition \eqref{eq:aci-step-assumptions} says that the increments of $1/\eta_t$ tend to zero.
The average of these increments then tends to zero as well, that is,
\[
  \frac{1}{t\,\eta_t}=\frac{1}{t\,\eta_1}+\frac1t\sum_{r<t}\Bigl(\frac{1}{\eta_{r+1}}-\frac{1}{\eta_r}\Bigr)\longrightarrow0 .
\]
Hence $t\eta_t\to\infty$, so $\sqrt{\eta_t}\ge t^{-1/2}$ for all large $t$, and $A^{\eta}_T\to\infty$.
With Step~2, this gives
\begin{equation}\label{eq:aci_mean_rate}
  0<\liminf_{T\to\infty}\frac{r_T(\calC^{\aci};\PP)}{A^{\eta}_T}
  \le\limsup_{T\to\infty}\frac{r_T(\calC^{\aci};\PP)}{A^{\eta}_T}<\infty,
\end{equation}
which is \eqref{eq:aci-sub}.
\end{proof}

\subsection{Proof of Corollary~\ref{cor:aci-decaying-risk}}
\label{app:aci-decaying-proof}

\begin{proof}
We check the step-size assumptions of \Cref{prop:aci-subopt} and then evaluate the sum in its conclusion.

\paragraph{Step 1: check the step-size assumptions.}
Let $c>0$ and $\beta\in(0,1)$.
The steps $\eta_t=ct^{-(1-\beta)}$ tend to zero because $\beta<1$, and
\[
 0\le\frac{1}{\eta_{t+1}}-\frac{1}{\eta_t}
 =\frac{(t+1)^{1-\beta}-t^{1-\beta}}{c}
 \le\frac{(1-\beta)t^{-\beta}}{c}\longrightarrow0
\]
because $\beta>0$.
Thus \Cref{prop:aci-subopt} applies.

\paragraph{Step 2: evaluate the cumulative rate.}
For $0\le\gamma<1$, comparison with an integral gives
\[
 T^{1-\gamma}\le\sum_{t=1}^Tt^{-\gamma}
 \le\int_0^T u^{-\gamma}\,du
 =\frac{T^{1-\gamma}}{1-\gamma}.
\]
Taking $\gamma=(1-\beta)/2$ yields
\[
 \sum_{t=1}^T\sqrt{\eta_t}
 =\sqrt c\sum_{t=1}^Tt^{-(1-\beta)/2}
 \asymp T^{(1+\beta)/2}.
\]
Substitution into \eqref{eq:aci-sub} proves \eqref{eq:aci-decaying-risk}.
\end{proof}

By \Cref{lem:aci_corridor}, these steps give the adversarial guarantee with $\kappa=(B+c)/c$.
For any prescribed $\kappa>1$, choose $c=B/(\kappa-1)$.
Since $(1+\beta)/2>\max\{1/2,\beta\}$ for $0<\beta<1$, the corollary gives a larger disagreement rate than truncated split conformal prediction in \Cref{thm:optimal} under the same guarantee.

\subsection{Proofs of Lemmas~\ref{lem:aci_risk} and~\ref{lem:aci_lip}}
\label{app:aci-simple-proofs}

\begin{proof}[Proof of \Cref{lem:aci_risk}]
The second equality in \eqref{eq:aci_risk} holds because $p_t-\alpha=-\chi(\Delta_t)$.
For the first equality, we prove two inequalities.

\emph{Lower bound.}
Let $\calC'\in\calO_\alpha(\PP)$ be any deterministic exact oracle.
Its oracle condition gives $\EE[\err_t^{\calC'}\mid\calG_{t-1}]=\alpha$ almost surely, so conditional Jensen's inequality gives
\[
 \EE_{\PP}\bigl[|\err_t-\err_t^{\calC'}|\,\big|\,\calG_{t-1}\bigr]
 \ge\bigl|\EE_{\PP}[\err_t-\err_t^{\calC'}\mid\calG_{t-1}]\bigr|=|p_t-\alpha| .
\]
Summing over $t$ and taking expectations, $\EE_{\PP}[\ell_T(\calC^{\aci},\calC')]\ge\EE_{\PP}\sum_{t\le T}|p_t-\alpha|$ for every oracle in $\calO_\alpha(\PP)$.

\emph{Upper bound.}
The population prediction set in \eqref{eq:population-oracle} is deterministic and measurable.
Since the fresh score is independent of the past and $F(q_{\PP,s})=1-\alpha$, it belongs to $\calO_\alpha(\PP)$.
The ACI set and the population set are thresholds of the same score, at the levels $q_t$ and $q_{\PP,s}$.
Their errors differ exactly when $s_t$ lies between the two levels, which has conditional probability $|F(q_t)-F(q_{\PP,s})|=|p_t-\alpha|$ given the past.
So this oracle attains the lower bound.
\end{proof}

\begin{proof}[Proof of \Cref{lem:aci_lip}]
By \Cref{lem:aci_range} and $q_{\PP,s}\in(0,B)$, we have $-\bar\eta\alpha-B<\Delta_t<B+\bar\eta(1-\alpha)$, so $|\Delta_t|\le B+\bar\eta$.
The function $\chi$ is nondecreasing with $\chi(0)=0$, so $\chi(u)\ge0$ for $u>0$ and $\chi(u)\le0$ for $u<0$.
It remains to prove the two bounds on $|\chi(u)|$.

For the upper bound, the extended distribution function $F$ has a density bounded by $\lambda_+$ on the whole real line, so $|\chi(u)|=|F(q_{\PP,s}+u)-F(q_{\PP,s})|\le\lambda_+|u|$ for every $u$.

For the lower bound, let $|u|\le B+\bar\eta$ and consider three cases.
If $q_{\PP,s}+u\in[0,B]$, then \eqref{eq:aci_lip} gives $|\chi(u)|\ge\lambda_-|u|\ge\underline{\lambda}|u|$.
If $q_{\PP,s}+u>B$, then $\chi(u)=1-(1-\alpha)=\alpha$, and $\underline{\lambda}|u|\le\nu\le\alpha$ because $|u|\le B+\bar\eta$.
If $q_{\PP,s}+u<0$, then $\chi(u)=-(1-\alpha)$, and $\underline{\lambda}|u|\le\nu\le1-\alpha$.
\end{proof}

\subsection{Proof of Lemma~\ref{lem:aci_fluctuation}}
\label{app:aci-fluctuation-proof}

The proof has two ingredients.
The first is a comparison lemma that turns a one-step inequality for a sequence into a bound on its size.
The second applies this lemma to the moments of $\Delta_t$.
The comparison is related to classical stochastic-approximation arguments \citep{chung1954stochastic,fabian1968asymptotic}.
We prove the form needed here for the general step sizes in \eqref{eq:aci-step-assumptions}.

\begin{lemma}[A comparison lemma for recursions]\label{lem:aci_rec}
Let $(\eta_t)$ be positive step sizes satisfying \eqref{eq:aci-step-assumptions}, let $m\ge1$ be an integer, and let $A>0$ and $D\ge0$.
Let $(z_t)$ be a nonnegative sequence.
Below, $o(\eta_t)$ and $o(1)$ stand for sequences of numbers whose ratio to $\eta_t$ (resp. whose value) tends to zero as $t\to\infty$.
\begin{enumerate}[label=\textnormal{(\alph*)}, leftmargin=2.2em, itemsep=0.2em]
  \item If $z_{t+1}\le\bigl(1-A\eta_t+o(\eta_t)\bigr)\,z_t+\bigl(D+o(1)\bigr)\,\eta_t^{m+1}$ for all large $t$, then $\limsup_{t\to\infty}z_t/\eta_t^m\le D/A$.
  \item If $z_{t+1}\ge\bigl(1-A\eta_t+o(\eta_t)\bigr)\,z_t+\bigl(D+o(1)\bigr)\,\eta_t^{m+1}$ for all large $t$, then $\liminf_{t\to\infty}z_t/\eta_t^m\ge D/A$.
\end{enumerate}
\end{lemma}

The level $(D/A)\eta_t^m$ is where contraction and added noise balance.
The two parts give the corresponding upper and lower comparisons.

\begin{lemma}[Moments of the deviation]\label{lem:aci_moment_bounds}
Suppose that the step sizes satisfy \eqref{eq:aci-step-assumptions}.
There exist constants $0<c_\Delta\le C_\Delta<\infty$ and a time $t_\star\ge1$ such that, for every $t\ge t_\star$,
\[
  c_\Delta\,\eta_t\le\EE\Delta_t^2\le C_\Delta\,\eta_t,
  \qquad
  \EE\Delta_t^4\le C_\Delta\,\eta_t^2.
\]
\end{lemma}

\begin{proof}[Proof of \Cref{lem:aci_fluctuation}]
We deduce the size of $\EE|\Delta_t|$ from the second and the fourth moment in \Cref{lem:aci_moment_bounds}.
Let $t\ge t_\star$.
For the upper bound, the Cauchy--Schwarz inequality gives
\[
  \EE|\Delta_t|\le(\EE\Delta_t^2)^{1/2}\le\sqrt{C_\Delta}\,\sqrt{\eta_t}.
\]
For the lower bound, we need that $\Delta_t$ is not concentrated on rare large values, and this is what the fourth moment controls.
H\"older's inequality gives $\EE\Delta_t^2\le(\EE|\Delta_t|)^{2/3}(\EE\Delta_t^4)^{1/3}$, and therefore
\[
  \EE|\Delta_t|
  \ge\frac{(\EE\Delta_t^2)^{3/2}}{(\EE\Delta_t^4)^{1/2}}
  \ge\frac{(c_\Delta\eta_t)^{3/2}}{(C_\Delta\eta_t^2)^{1/2}}
  =\frac{c_\Delta^{3/2}}{C_\Delta^{1/2}}\,\sqrt{\eta_t}.
\]
This proves the lemma with $C_0=C_\Delta^{1/2}$ and $c_0=c_\Delta^{3/2}/C_\Delta^{1/2}$.
\end{proof}

\subsection{Proof of Lemma~\ref{lem:aci_rec}}
\label{app:aci-recursion-proof}

\begin{proof}
We first divide the sequence by $\eta_t^m$, and then compare the result with the level at which the two effects balance.

\paragraph{Step 1: normalization.}
Let $\tilde z_t=z_t/\eta_t^m$.
By \eqref{eq:aci-step-assumptions},
\[
  \frac{\eta_t}{\eta_{t+1}}=1+\eta_t\Bigl(\frac{1}{\eta_{t+1}}-\frac{1}{\eta_t}\Bigr)=1+o(\eta_t),
\]
and, since the step sizes are bounded, also $(\eta_t/\eta_{t+1})^m=1+o(\eta_t)$.
So the change of the normalization from $\eta_t^m$ to $\eta_{t+1}^m$ only contributes a factor of this form.
Dividing the inequality in (a) by $\eta_{t+1}^m$ and multiplying out, we get
\begin{equation}\label{eq:aci_rec_normalized}
  \tilde z_{t+1}\le\bigl(1-A\eta_t+o(\eta_t)\bigr)\,\tilde z_t+\bigl(D+o(1)\bigr)\,\eta_t
  \qquad\text{for all large }t .
\end{equation}
Only the deterministic coefficients are expanded here; no bound on $\tilde z_t$ is assumed.
In the same way, the inequality in (b) gives $\tilde z_{t+1}\ge\bigl(1-A\eta_t+o(\eta_t)\bigr)\,\tilde z_t+\bigl(D+o(1)\bigr)\,\eta_t$ for all large $t$.

\paragraph{Step 2: the step sizes are not summable.}
The signed increments of $1/\eta_t$ tend to zero, so their Ces\`aro averages give
\[
 \frac1{t\eta_t}=\frac1{t\eta_1}
 +\frac1t\sum_{r<t}(\eta_{r+1}^{-1}-\eta_r^{-1})\longrightarrow0.
\]
Thus $t\eta_t\to\infty$, and $\sum_t\eta_t=\infty$.

\paragraph{Step 3: part (a).}
Fix $\varepsilon\in(0,A)$.
For all large $t$, the two $o$-terms in \eqref{eq:aci_rec_normalized} are at most $\varepsilon\eta_t$ and $\varepsilon$.
Since $\tilde z_t\ge0$, this gives, for all large $t$,
\[
  \tilde z_{t+1}\le\bigl(1-(A-\varepsilon)\eta_t\bigr)\,\tilde z_t+(D+\varepsilon)\,\eta_t .
\]
Let $L^+=(D+\varepsilon)/(A-\varepsilon)$ be the level at which the two terms balance.
Subtracting $L^+$ from both sides gives $\tilde z_{t+1}-L^+\le(1-(A-\varepsilon)\eta_t)(\tilde z_t-L^+)$.
The factor lies in $(0,1)$ for large $t$, so the positive part satisfies
\[
  (\tilde z_{t+1}-L^+)_+\le\bigl(1-(A-\varepsilon)\eta_t\bigr)(\tilde z_t-L^+)_+ .
\]
Iterating from a large time $t_1$,
\[
  (\tilde z_T-L^+)_+\le(\tilde z_{t_1}-L^+)_+\prod_{t=t_1}^{T-1}\bigl(1-(A-\varepsilon)\eta_t\bigr)\longrightarrow0
  \qquad\text{as }T\to\infty ,
\]
because $\sum_t\eta_t=\infty$.
Hence $\limsup_t\tilde z_t\le L^+$, and letting $\varepsilon\downarrow0$ gives $\limsup_t\tilde z_t\le D/A$.

\paragraph{Step 4: part (b).}
$D=0$ directly implies $\tilde z_t\ge0$.
Otherwise fix $\varepsilon\in(0,D)$ and let $L^-=(D-\varepsilon)/(A+\varepsilon)$.
For all large $t$, the two $o$-terms are at least $-\varepsilon\eta_t$ and $-\varepsilon$, so the lower inequality gives $\tilde z_{t+1}\ge(1-(A+\varepsilon)\eta_t)\,\tilde z_t+(D-\varepsilon)\,\eta_t$.
As in Step~3, this gives $(L^--\tilde z_{t+1})_+\le(1-(A+\varepsilon)\eta_t)(L^--\tilde z_t)_+$ for all large $t$.
Hence $(L^--\tilde z_T)_+\to0$, so $\liminf_t\tilde z_t\ge L^-$, and letting $\varepsilon\downarrow0$ gives $\liminf_t\tilde z_t\ge D/A$.
\end{proof}

\subsection{Proof of Lemma~\ref{lem:aci_moment_bounds}}
\label{app:aci-moments-proof}

\begin{proof}
The idea is to raise the recursion \eqref{eq:aci_delta} to the second and the fourth power and to take expectations.
The noise has mean zero given the past, so the terms that are linear in the noise vanish.
What remains is a one-step inequality for each moment, to which \Cref{lem:aci_rec} applies.

We write $\Delta_{t+1/2}=\Delta_t-\eta_t\chi(\Delta_t)$ for the deviation, so that $\Delta_{t+1}=\Delta_{t+1/2}+\eta_t\xi_t$.
The variable $\Delta_{t+1/2}$ is determined by the past.
The range bound and $\bar\eta<\infty$ make $\Delta_t$, $\Delta_{t+1/2}$, and all the moments below bounded.

\paragraph{Step 1: contraction.}
By \Cref{lem:aci_lip}, $\underline{\lambda}\Delta_t^2\le\Delta_t\chi(\Delta_t)\le\lambda_+\Delta_t^2$ and $\chi(\Delta_t)^2\le\lambda_+^2\Delta_t^2$.
Expanding $\Delta_{t+1/2}^2=\Delta_t^2-2\eta_t\Delta_t\chi(\Delta_t)+\eta_t^2\chi(\Delta_t)^2$ therefore gives
\begin{equation}\label{eq:aci_Ysq}
  (1-2\lambda_+\eta_t)\,\Delta_t^2
  \;\le\; \Delta_{t+1/2}^2
  \;\le\;(1-2\underline{\lambda}\eta_t+\lambda_+^2\eta_t^2)\,\Delta_t^2 .
\end{equation}

\paragraph{Step 2: the second moment, upper bound.}
Squaring $\Delta_{t+1}=\Delta_{t+1/2}+\eta_t\xi_t$ and using \eqref{eq:aci_noise},
\begin{equation}\label{eq:aci_second_cond}
  \EE\bigl[\Delta_{t+1}^2\mid\calG_{t-1}\bigr]=\Delta_{t+1/2}^2+\eta_t^2\,p_t(1-p_t).
\end{equation}
Taking expectations and using the upper bound in \eqref{eq:aci_Ysq} and $p_t(1-p_t)\le1/4$,
\[
  \EE\Delta_{t+1}^2\le(1-2\underline{\lambda}\eta_t+\lambda_+^2\eta_t^2)\,\EE\Delta_t^2+\tfrac14\eta_t^2 .
\]
\Cref{lem:aci_rec}(a) with $z_t=\EE\Delta_t^2$, $m=1$, $A=2\underline{\lambda}$, and $D=1/4$ gives $\limsup_t\EE\Delta_t^2/\eta_t\le1/(8\underline{\lambda})$.
In particular $\EE\Delta_t^2\to0$.

\paragraph{Step 3: the second moment, lower bound.}
For the lower bound we need that the variance of the noise does not vanish.
For $p\in[0,1]$ we have $|p(1-p)-\alpha(1-\alpha)|=|p-\alpha|\,|1-p-\alpha|\le|p-\alpha|$, and $|p_t-\alpha|=|\chi(\Delta_t)|\le\lambda_+|\Delta_t|$.
Hence
\[
  \bigl|\EE[p_t(1-p_t)]-\sigma^2\bigr|\le\lambda_+\EE|\Delta_t|\le\lambda_+\bigl(\EE\Delta_t^2\bigr)^{1/2}\longrightarrow0
\]
by Step~2.
Taking expectations in \eqref{eq:aci_second_cond} and using the lower bound in \eqref{eq:aci_Ysq},
\[
  \EE\Delta_{t+1}^2\ge(1-2\lambda_+\eta_t)\,\EE\Delta_t^2+\bigl(\sigma^2+o(1)\bigr)\,\eta_t^2 .
\]
\Cref{lem:aci_rec}(b) with $z_t=\EE\Delta_t^2$, $m=1$, $A=2\lambda_+$, and $D=\sigma^2$ gives $\liminf_t\EE\Delta_t^2/\eta_t\ge\sigma^2/(2\lambda_+)>0$.
Steps~2 and~3 prove the two bounds on $\EE\Delta_t^2$.

\paragraph{Step 4: the fourth moment.}
We expand $(\Delta_{t+1/2}+\eta_t\xi_t)^4$ and take the conditional expectation given the past.
The term that is linear in $\xi_t$ vanishes.
For the other terms we use $\EE[|\xi_t|^j\mid\calG_{t-1}]\le\EE[\xi_t^2\mid\calG_{t-1}]\le1/4$ for $j\ge2$, which holds because $|\xi_t|\le1$.
This gives
\[
  \EE\bigl[\Delta_{t+1}^4\mid\calG_{t-1}\bigr]
  \le \Delta_{t+1/2}^4+\tfrac32\eta_t^2\Delta_{t+1/2}^2+\eta_t^3|\Delta_{t+1/2}|+\tfrac14\eta_t^4 .
\]
Squaring the upper bound in \eqref{eq:aci_Ysq} gives $\Delta_{t+1/2}^4\le(1-4\underline{\lambda}\eta_t+C\eta_t^2)\Delta_t^4$.
Moreover, $\EE \Delta_{t+1/2}^2\le C\eta_t$ for all large $t$ by Step~2, and $|\Delta_{t+1/2}|$ is bounded.
So the last three terms have expectation at most $C\eta_t^3$, and for all large $t$,
\[
  \EE\Delta_{t+1}^4\le(1-4\underline{\lambda}\eta_t+C\eta_t^2)\,\EE\Delta_t^4+C\eta_t^3 .
\]
\Cref{lem:aci_rec}(a) with $z_t=\EE\Delta_t^4$, $m=2$, $A=4\underline{\lambda}$, and $D=C$ gives $\limsup_t\EE\Delta_t^4/\eta_t^2<\infty$.

The three limits in Steps~2--4 give the lemma, with suitable constants $c_\Delta$ and $C_\Delta$ and a suitable time $t_\star$.
\end{proof}

\end{document}